\documentclass[a4paper]{article}
\usepackage[inner=30mm,outer=30mm,textheight=225mm]{geometry}
\usepackage{amsmath}
\usepackage{amssymb}
\usepackage{array}
\usepackage{footnote}
\usepackage{graphicx}
\usepackage{psfrag}
\usepackage{theorem}

\theorembodyfont{\upshape}
\newtheorem{theorem}{Theorem}[section]

\newtheorem{defn}[theorem]{Definition}

\newtheorem{example}[theorem]{Example}
\newtheorem{lemma}[theorem]{Lemma}

\newtheorem{algorithm}[theorem]{Algorithm}

\newcommand{\qed}{\hfill\rule[-0.5mm]{1.5mm}{3.0mm}}
\newtheorem{proofthm}{Proof}
\newenvironment{proof}{\begin{proofthm}}{\qed \end{proofthm}}

\newcommand{\discref}{\bar{D}}
\renewcommand{\epsilon}{\varepsilon}

\newcommand{\kb}{K^2}
\newcommand{\lst}{\mathrm{LST}}

\newcommand{\mobius}{M\"{o}bius}

\newcommand{\ppirr}{{$\mathbb{P}^2$-ir\-re\-du\-ci\-ble}}
\newcommand{\ppirrty}{{$\mathbb{P}^2$-ir\-re\-du\-ci\-bi\-li\-ty}}

\newcommand{\R}{\mathbb{R}}
\newcommand{\regina}{{\em Regina}}

\newcommand{\rpp}{\mathbb{R}P^2}
\newcommand{\rps}{\mathbb{R}P^3}

\newcommand{\sfs}[2]{\mathrm{SFS}\left(#1: #2\right)}
\newcommand{\sfslong}{Seifert fibred space}

\newcommand{\scircle}{S^1}

\newcommand{\torus}{T^2}

\newcommand{\twt}{\tilde{T}}

\newcommand{\Z}{\mathbb{Z}}

\newcommand{\homtwo}[4]{
    \mbox{\scriptsize \renewcommand{\arraystretch}{1}
        $\! \left[ \begin{array}{@{\ }c@{\ }c@{\ }} #1 & #2 \\ #3 & #4
        \end{array} \right]$
    }
}

\newcommand{\homtwolarge}[4]{
    \left[ \begin{array}{cc} #1 & #2 \\ #3 & #4
    \end{array} \right]
}

\title{Structures of small closed non-orientable 3-manifold triangulations}
\author{Benjamin A.~Burton}
\date{September 15, 2005\footnote{The initial version of this paper was
released in November 2003.  The update from September 2005 corrects an
off-by-one error in the formulae of
Theorems~\ref{t-idpluggedibundletorus} and~\ref{t-idpluggedthick}.
It also includes general proofreading, particularly in the discussion
of layered solid tori.}}

\begin{document}

\maketitle

\abstract{A census is presented of all closed non-orientable 3-manifold
triangulations formed from at most seven tetrahedra satisfying the
additional constraints of minimality and {\ppirrty}.  The eight
different 3-manifolds represented by these 41 different triangulations
are identified and described in detail, with particular attention paid
to the recurring combinatorial structures that are shared amongst the
different triangulations.  Using these recurring structures,
the resulting triangulations are generalised to infinite families that
allow similar triangulations of additional 3-manifolds to be formed.
Algorithms and techniques used in constructing the census are included.}

\section{Introduction} \label{s-intro}

It is useful when studying 3-manifold topology to have a complete reference
of all 3-manifold triangulations satisfying some broad set of constraints.
Examples include Callahan, Hildebrand and Weeks'
census of cusped hyperbolic 3-manifold triangulations formed from at most
seven tetrahedra \cite{cuspedcensus} and Matveev's census of closed
orientable triangulations formed from at most six tetrahedra \cite{matveev6}.

Such references provide an excellent pool of examples for testing hypotheses
and searching for triangulations that satisfy unusual properties.  In
addition they offer insight into the structures of minimal 3-manifold
triangulations.  In fact, since very few sufficient conditions for minimality
are currently known, censuses play an important role in proving the minimality
of small triangulations.

Much recent progress has been made in enumerating closed orientable
3-manifolds and their triangulations.  Matveev presents
a census of closed orientable triangulations formed from at most six
tetrahedra \cite{matveev6},
extended to seven tetrahedra by Ovchinnikov.  Martelli and Petronio
form a census of closed orientable 3-manifolds formed
from up to nine tetrahedra \cite{italian9},
although they are primarily concerned with the
3-manifolds and their geometric structures and so the triangulations
themselves are not listed.  More recently their census
has been extended to ten tetrahedra by Martelli \cite{italian10}.

Less progress has been made regarding closed non-orientable
triangulations.  Amendola and Martelli present a
census of closed non-orientable 3-manifolds formed from up to six
tetrahedra \cite{italian-nor6}, a particularly interesting census because
it is constructed
without the assistance of a computer.  Again these authors are primarily
concerned with the 3-manifolds and their geometric structures, and so
not all triangulations of these 3-manifolds are obtained.

Here we extend the closed non-orientable census of Amendola and Martelli
to seven tetrahedra, and in addition we enumerate the many different
triangulations of these 3-manifolds instead of just the 3-manifolds
themselves.
Furthermore we examine the combinatorial structures of these
triangulations in detail, highlighting common constructions
that recur throughout the census triangulations.

Independently of this work, Amendola and Martelli
have also announced an extension of their theoretical census to seven
tetrahedra \cite{italian-nor7}.  As in their six-tetrahedron census they
concentrate only on the resulting 3-manifolds and not their different
triangulations, but again the non-computational nature of their work is
remarkable.

As with the previous closed censuses described above,
we consider only triangulations satisfying the following constraints.
\begin{itemize}
    \item {\em Closed:}  The triangulation is of a closed
    3-manifold.  In particular it has no boundary faces, and each vertex
    link is a 2-sphere.
    \item {\em \ppirr:}  The underlying 3-manifold has no
    embedded two-sided projective planes, and furthermore every
    embedded 2-sphere bounds a ball.
    \item {\em Minimal:}  The underlying 3-manifold cannot be
    triangulated using strictly fewer tetrahedra.
\end{itemize}

Requiring triangulations to be {\ppirr} and minimal keeps the number of
triangulations down to manageable levels, focussing only upon the simplest
triangulations of the simplest 3-manifolds (from which more complex
3-manifolds can be constructed).
Minimal triangulations prove to be particularly useful for studying
the 3-manifolds that they represent, since they are
frequently well structured as seen in both Matveev's census
\cite{matveev6} and the results presented here.

For the $\leq 7$-tetrahedron non-orientable census presented in this paper,
a brief summary of results is presented in Table~\ref{tab-summary}.
Each triangulation is counted once up to isomorphism, i.e., a
relabelling of the tetrahedra within the triangulation
and their individual faces.  It is worth
noting that the number of triangulations is significantly larger than
the number of 3-manifolds, since most 3-manifolds in the census can
be realised by several different minimal triangulations.

\begin{table}[htb]
\begin{center} \begin{tabular}{|r|r|r|}
    \hline
    \bf Tetrahedra & \bf 3-Manifolds & \bf Triangulations \\
    \hline
    $\leq 5$ & 0 & 0 \\
    6 & 5 & 24 \\
    7 & 3 & 17 \\
    \hline
    Total & 8 & 41 \\
    \hline
\end{tabular} \end{center}
\caption{Summary of closed non-orientable census results}
\label{tab-summary}
\end{table}

The final 41 triangulations are found to be remarkably similar in their
construction.  By identifying these similarities we construct a small
handful of infinite parameterised families of 3-manifold triangulations
that encompass 38 of these 41 triangulations.  The remaining three
triangulations are all six-tetrahedron triangulations and might well be
small exceptional cases that do not generalise at all --- an extension
of this census to higher numbers of tetrahedra should offer further
insight.

In Section~\ref{s-algm} we describe the method by which this census was
constructed.  The remainder of this paper is devoted to presenting the
census results and describing in detail the combinatorial structures of
the various triangulations.  Section~\ref{s-common} describes the
construction of thin $I$-bundles and layered solid tori, which are
parameterised building blocks that recur frequently throughout the census
triangulations.  In Section~\ref{s-families} we combine these building
blocks to form our infinite families of 3-manifold
triangulations, as well as describing the three exceptional
triangulations from the census that these families do not cover.
Finally Section~\ref{s-census} closes with a full
listing of the 41 triangulations found in the census, using the
constructions of Sections~\ref{s-common} and~\ref{s-families} to
simplify their descriptions and identify the underlying 3-manifolds.

All of the computational work was performed using {\regina},
a computer program that performs a variety of different calculations and
procedures in 3-manifold topology \cite{regina,burton-regina}.
The program \regina, its source code
and accompanying documentation are freely available from
{\tt http://\allowbreak regina.\allowbreak sourceforge.\allowbreak net/}.

Special thanks must go to J.~Hyam Rubinstein for many helpful
discussions throughout the course of this research.
The author would also like to thank
the Australian Research Council, RMIT University and
the University of Melbourne for their support.

\section{Constructing the Census} \label{s-algm}

As with most of the prior censuses listed in Section~\ref{s-intro}, this
census is based upon a computer search through the possible
triangulations that can be formed from various numbers of tetrahedra.
At the heart of this computer search is a procedure that generates
triangulations from $n$ tetrahedra using all possible identifications of
tetrahedron faces under all possible rotations and reflections.  Note
that in practice this search can be refined using combinatorial
techniques to avoid many isomorphic duplicates of triangulations.

Recall from Section~\ref{s-intro} that we require only triangulations
that are closed, non-orientable, minimal and {\ppirr}.  Some of these
properties are easily incorporated into the generation procedure
described above.  For instance, it is straightforward to adjust the
generation procedure so that only non-orientable triangulations with no
boundary faces are produced (we simply ensure that
every tetrahedron face is matched with a partner, and we
keep track of the orientability of each tetrahedron
as we go).  Other properties such as minimality and {\ppirrty} are
less easily dealt with, and in many cases must be evaluated
after potential triangulations have already been constructed.

Thus we can decompose the construction of a census into the two stages of
generation and analysis as follows.

\begin{algorithm}[Census Construction] \label{a-census}
    Let $n$ be some positive integer.
    A census of all closed non-orientable minimal {\ppirr}
    triangulations formed from $n$ tetrahedra can be constructed in the
    following fashion.
    \begin{enumerate}
        \item \label{en-algmgeneration}
        Generate a set of triangulations that is guaranteed to include
        everything that should appear in the census.
        For instance, we might generate all closed
        non-orientable 3-manifold triangulations formed from $n$
        tetrahedra without regard for minimality or {\ppirrty}.
        In practise the set of triangulations that we generate
        is more complicated,
        as seen in Section~\ref{s-algmgeneration} below.

        It should be ensured at this stage that no two triangulations in the
        generated set are isomorphic duplicates.

        \item \label{en-algmanalysis}
        For each triangulation generated in step~\ref{en-algmgeneration},
        determine whether it satisfies the full set of census constraints
        and if so then include it in the final results.  For instance, in
        the example above we would need to test each triangulation for
        minimality and {\ppirrty}.
    \end{enumerate}
\end{algorithm}

We proceed to examine in detail the individual steps of the census algorithm
in Sections~\ref{s-algmgeneration} and~\ref{s-algmanalysis} below.

\subsection{Generation of Triangulations} \label{s-algmgeneration}

It can be observed that a triangulation formed from $n$ tetrahedra can be
uniquely determined by the following information:
\begin{itemize}
    \item A {\em face pairing}, i.e., a partition of the $4n$
    tetrahedron faces into $2n$ pairs indicating which tetrahedron faces
    are to be identified with which others;
    \item A list of the $2n$ rotations and/or reflections used to
    perform each of these $2n$ face identifications.
\end{itemize}

A first approach to step~\ref{en-algmgeneration} of
Algorithm~\ref{a-census}, i.e., the generation of triangulations,
could thus be as follows.
\begin{enumerate}
    \item \label{en-algmgenpairings}
    Enumerate all possible face pairings.

    \item \label{en-algmgengluings}
    For each face pairing, try all possible $6^{2n}$ combinations
    of rotations and reflections for the $2n$ face identifications.
    Record each triangulation thus produced.
\end{enumerate}

Note that there is no guarantee that the triangulations produced are
actually triangulations of 3-manifolds.  They are however guaranteed to have
no boundary faces, and by keeping track of the orientation of each
tetrahedron as we go it is straightforward to ensure non-orientability.

In practice this generation is exceptionally slow; we see already that
step~\ref{en-algmgengluings} is exponential in the number of tetrahedra.
This accounts for the limited extent of current census results described in
Section~\ref{s-intro}.  By exploiting the fact that we are interested
only in closed non-orientable minimal {\ppirr} triangulations of
3-manifolds, we can use a variety of methods to improve the running time
of the census algorithm.  These methods include the following.
\begin{itemize}
    \item Several results relating to face
    pairings are presented in \cite{burton-facegraphs}, which allow many
    of the face pairings to be tossed away in
    step~\ref{en-algmgenpairings} and which use properties of the remaining
    face pairings to provide significant improvements to
    step~\ref{en-algmgengluings}.  For the six-tetrahedron non-orientable
    census these results eliminate over 98\% of the running time of the
    census generation.

    \item Further improvements can be made by modifying
    step~\ref{en-algmgengluings} to avoid low degree edges, a technique
    used successfully in earlier hyperbolic censuses
    \cite{cuspedcensus,cuspedcensusold} and closed
    orientable censuses \cite{italian9,matveev6}.
    Details of how this technique is applied to a closed non-orientable
    census can be found in \cite{burton-facegraphs}.
\end{itemize}

The issue of isomorphism must also be dealt with.
Isomorphic duplicates of face pairings are avoided by selecting a
canonical representation for each face pairing from amongst all possible
relabellings.  For instance, the canonical representation might be the
lexicographically smallest relabelling when written in some standard format.

Any face pairing that is not in its canonical representation can therefore
be ignored.  Furthermore, the generation of face pairings can be
streamlined to toss away partially constructed face pairings that will
clearly not have this property.  Isomorphic duplicates of entire triangulations
are dealt with in a similar way.

\subsection{Analysis of Triangulations} \label{s-algmanalysis}

For each triangulation $T$ that is constructed in
step~\ref{en-algmgeneration} of Algorithm~\ref{a-census}, we must still
determine whether $T$ is in fact a closed non-orientable minimal {\ppirr}
triangulation of a 3-manifold.
It is straightforward to test whether $T$ is indeed a
3-manifold triangulation, and the generation described in
Section~\ref{s-algmgeneration} is already tweaked to ensure
closedness and non-orientability.

It remains then to determine whether each triangulation is
minimal and {\ppirr}.  These properties are somewhat more difficult
to test for.  No general test for minimality is currently known.
Furthermore, current tests for reducibility involve either cutting along
embedded 2-spheres (see for instance \cite{jacotollefson-algorithms}) or
crushing embedded 2-spheres to a
point (see \cite{jacorubin-algorithms}).  Cutting along 2-spheres produces
very large triangulations that make such algorithms too slow for
practice, and it is difficult to generalise the crushing
algorithms to the non-orientable case due to a number of additional
complications that arise.

We can however use a variety of alternative techniques to at least
answer the questions of minimality and {\ppirrty} for some
triangulations.  These techniques are outlined in
Sections~\ref{s-algmanalysis-eltmoves} to~\ref{s-algmanalysis-inv}
below.  Section~\ref{s-algmanalysis-results} concludes with a
discussion of the success of these techniques when applied to the particular
$\leq 7$-tetrahedron non-orientable census under consideration.

\subsubsection{Elementary Moves} \label{s-algmanalysis-eltmoves}

It is often possible to make a modification local to a few tetrahedra
within a triangulation that preserves the underlying 3-manifold, with no
knowledge whatsoever of the global triangulation structure or of any
properties of the 3-manifold.  Such local modifications are referred to
as {\em elementary moves}.

If a sequence of elementary moves can be found that reduces the
number of tetrahedra in a triangulation, then since these moves
preserve the underlying 3-manifold it follows that the original
triangulation cannot be minimal.

Alternatively, if a sequence of elementary moves can be found that
links two census triangulations with the same number of tetrahedra,
it follows that any deductions regarding the minimality or {\ppirrty}
of one triangulation can be simultaneously applied to the other.

The individual elementary moves that were used for this census
are described below.  These moves are
not new, and many were implemented in 1999 by David Letscher in his
computer program {\em Normal}.  Analogous techniques using local
modifications of special spines have been used by Matveev
\cite{matveev6} and Martelli and Petronio \cite{italian9}.

Note that for each of these elementary moves, it is
relatively straightforward to test whether the move may be made.  No
properties are assumed of the underlying triangulation, and so these
moves may be safely applied to triangulations with boundary faces or
ideal vertices (vertices with higher genus links).

\begin{lemma}[Pachner Moves]
    If a triangulation contains a non-boundary edge of degree three
    that belongs to three distinct tetrahedra then the following move
    can be made.

    These three tetrahedra (adjacent along our edge of degree three) are
    replaced with a pair of tetrahedra adjacent along a single face,
    as illustrated in Figure~\ref{fig-eltmove32}.
    This is called a {\em 3-2 Pachner move} and is one of the bistellar
    operations considered in \cite{pachner-moves}.
    This move preserves the underlying 3-manifold and reduces the number
    of tetrahedra in its triangulation by one.

    \begin{figure}[htb]
    \centerline{\includegraphics[scale=0.7]{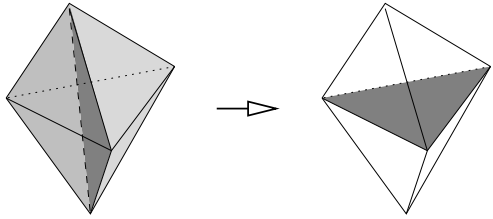}}
    \caption{A 3-2 Pachner move}
    \label{fig-eltmove32}
    \end{figure}

    Likewise, if a triangulation contains a non-boundary face belonging
    to two distinct tetrahedra then this move can be applied in reverse.
    This is called a {\em 2-3 Pachner move} and again preserves the
    underlying 3-manifold, this time increasing the number of tetrahedra
    by one.
\end{lemma}

\begin{proof}
    Since the modifications for each move take place entirely
    within the boundary of the
    illustrated polyhedron and since none of the vertex or edge links on
    the boundary are changed, it is clear that the underlying 3-manifold
    is preserved.
\end{proof}

\begin{lemma}[4-4 Move]
    If a triangulation contains a non-boundary edge of degree four
    that belongs to four distinct tetrahedra then the following move
    can be made.

    These four tetrahedra are replaced with four different tetrahedra
    meeting along a new edge of degree four, running perpendicular to
    the old edge of degree four.  This is called a {\em 4-4 move} and is
    illustrated in Figure~\ref{fig-eltmove44}.  This move preserves both
    the underlying 3-manifold and the number of tetrahedra.

    \begin{figure}[htb]
    \centerline{\includegraphics[scale=0.7]{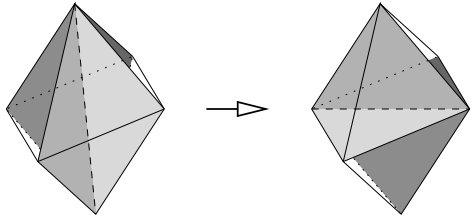}}
    \caption{A 4-4 move}
    \label{fig-eltmove44}
    \end{figure}
\end{lemma}

\begin{proof}
    Again the modifications take place entirely within the illustrated
    polyhedron, with no changes to the boundary vertex and edge links.
    Therefore the underlying 3-manifold is preserved.  Note that a 4-4
    move is simply a 2-3 Pachner move followed by a 3-2 Pachner move.
\end{proof}

\begin{lemma}[2-0 Vertex Move] \label{l-eltmove20v}
    Let $v$ be an internal vertex of degree two in a triangulation.
    Assume that the two tetrahedra meeting
    $v$ are distinct, and that these two tetrahedra
    meet along three different faces as
    illustrated in the left hand diagram of Figure~\ref{fig-eltmove20v}.
    Assume that the remaining faces of each tetrahedron, i.e., the
    two faces opposite $v$, are distinct
    and are not both boundary faces (though one of them may be a
    boundary face).

    \begin{figure}[htb]
    \psfrag{V}{{\small $v$}}
    \centerline{\includegraphics[scale=0.7]{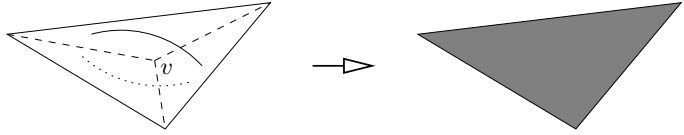}}
    \caption{A 2-0 vertex move}
    \label{fig-eltmove20v}
    \end{figure}

    Then these two tetrahedra may be flattened to a single face as
    illustrated in the right hand diagram of
    Figure~\ref{fig-eltmove20v}.  This is called a {\em 2-0 vertex move}.
    This move preserves the underlying 3-manifold and reduces the number
    of tetrahedra in its triangulation by two.
\end{lemma}

\begin{proof}
    Since the two faces opposite $v$ are distinct and are
    not both boundary faces, there is some third tetrahedron
    $\Delta$ adjacent to one of these faces as illustrated in the
    left hand diagram of Figure~\ref{fig-eltmove20vexpl}.

    \begin{figure}[htb]
    \psfrag{D}{$\Delta$}
    \psfrag{V}{{\small $v$}}
    \centerline{\includegraphics[scale=0.7]{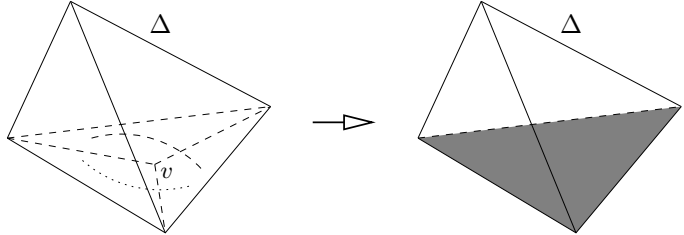}}
    \caption{A 2-0 vertex move explained}
    \label{fig-eltmove20vexpl}
    \end{figure}

    The result of performing the 2-0 vertex move
    is illustrated in the right hand diagram of
    Figure~\ref{fig-eltmove20vexpl}.  The modifications all take place
    within the boundary of the illustrated polyhedron
    and none of the vertex or edge links on the boundary of this
    polyhedron are changed.  Thus the 2-0 vertex move preserves the
    underlying 3-manifold.
\end{proof}

\begin{lemma}[2-0 Edge Move] \label{l-eltmove20e}
    Let $e$ be a non-boundary edge of degree two in a triangulation, as
    illustrated in the left hand diagram of Figure~\ref{fig-eltmove20e}.
    Assume that the two tetrahedra meeting $e$ are distinct.  Assume also
    that the edges opposite $e$ in each tetrahedron, labelled $g$ and
    $h$ in the diagram, are distinct and are not both boundary edges
    (though one of them may be a boundary edge).

    \begin{figure}[htb]
    \psfrag{G1}{{\small $G_1$}} \psfrag{G2}{{\small $G_2$}}
    \psfrag{H1}{{\small $H_1$}} \psfrag{H2}{{\small $H_2$}}
    \psfrag{e}{{\small $e$}} \psfrag{g}{{\small $g$}} \psfrag{h}{{\small $h$}}
    \centerline{\includegraphics[scale=0.7]{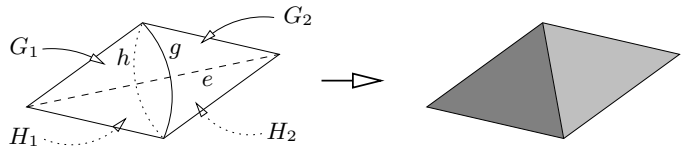}}
    \caption{A 2-0 edge move}
    \label{fig-eltmove20e}
    \end{figure}

    Consider now the four faces in the diagram that do not contain edge
    $e$ (these are the faces on either side of edges $g$ and $h$).  Label
    these faces $G_1$, $G_2$, $H_1$ and $H_2$ as illustrated, so that
    edge $g$ lies between faces $G_1$ and $G_2$ and edge $h$ lies
    between faces $H_1$ and $H_2$.
    Assume that faces $G_1$ and $H_1$ are distinct and that faces
    $G_2$ and $H_2$ are distinct.  Assume that
    all four of these faces are not identified in pairs (though we may
    have two of these faces identified, such as $G_1$ and $G_2$).
    Assume that we do not have a situation in which two of these faces
    are identified with each other and the remaining two are boundary faces.

    Then these two tetrahedra may be flattened to a pair of faces, as
    illustrated in the right hand diagram of Figure~\ref{fig-eltmove20e}.
    This is called a {\em 2-0 edge move}.  This move preserves the
    underlying 3-manifold and reduces the number of tetrahedra in its
    triangulation by two.
\end{lemma}

\begin{proof}
    Consider the disc bounded by edges $g$ and $h$ that slices through
    our two tetrahedra.  Since $g$ and $h$ are distinct and are not
    both boundary edges, we may crush this disc to a single edge without
    changing the underlying 3-manifold.
    After this operation we are left with two triangular pillows
    joined along a single edge as
    illustrated in the bottom left hand diagram of
    Figure~\ref{fig-eltmove20eexpl}.
    We may then retriangulate each of these pillows using two
    tetrahedra as illustrated in the bottom right hand diagram of
    Figure~\ref{fig-eltmove20eexpl}.

    \begin{figure}[htb]
    \psfrag{G1}{{\small $G_1$}} \psfrag{G2}{{\small $G_2$}}
    \psfrag{H1}{{\small $H_1$}} \psfrag{H2}{{\small $H_2$}}
    \psfrag{e}{{\small $e$}} \psfrag{g}{{\small $g$}} \psfrag{h}{{\small $h$}}
    \centerline{\includegraphics[scale=0.7]{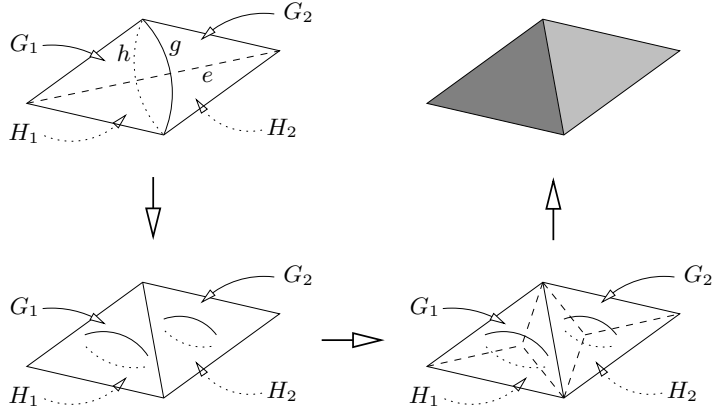}}
    \caption{The intermediate stages of a 2-0 edge move}
    \label{fig-eltmove20eexpl}
    \end{figure}

    The conditions placed upon faces $G_1$, $G_2$, $H_1$ and $H_2$
    allow us to use Lemma~\ref{l-eltmove20v} to perform a
    2-0 vertex move upon each pillow, thereby flattening each pillow to
    a face.  This procedure is illustrated in the top right hand diagram of
    Figure~\ref{fig-eltmove20eexpl}.  This completes the 2-0 edge
    move, leaving no changes in the underlying 3-manifold and an overall
    reduction of two tetrahedra.
\end{proof}

\begin{lemma}[2-1 Edge Move]
    Let $e$ be a non-boundary edge of degree one in a triangulation, as
    illustrated in the left hand diagram of Figure~\ref{fig-eltmove21e}.
    Label the single tetrahedron containing $e$ as $\Delta$, label the
    endpoints of $e$ as $A$ and $B$ and label the remaining two vertices
    of $\Delta$ as $C$ and $D$.

    \begin{figure}[htb]
    \psfrag{A}{{\small $A$}} \psfrag{B}{{\small $B$}}
    \psfrag{C}{{\small $C$}} \psfrag{D}{{\small $D$}}
    \psfrag{E}{{\small $E$}}
    \psfrag{e}{{\small $e$}}
    \psfrag{g}{{\small $g$}} \psfrag{h}{{\small $h$}}
    \psfrag{T}{{\small $\Delta$}} \psfrag{S}{{\small $\Delta'$}}
    \centerline{\includegraphics[scale=0.7]{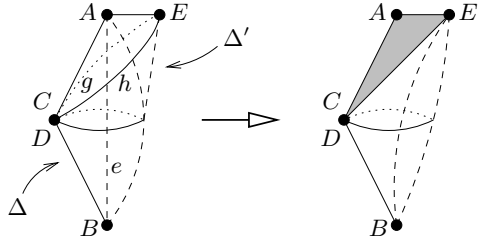}}
    \caption{A 2-1 edge move}
    \label{fig-eltmove21e}
    \end{figure}

    Assume that face {\em CAD} (the upper face of $\Delta$)
    is not a boundary face and let $\Delta'$ be the
    tetrahedron adjacent along this face.  Assume that
    $\Delta$ and $\Delta'$ are distinct tetrahedra and label the
    remaining vertex of $\Delta'$ as $E$.
    Assume that edges {\em CE} and {\em DE} of tetrahedron $\Delta'$
    are distinct and label them $g$ and $h$ respectively.
    Assume that edges $g$ and $h$ are not both boundary edges
    (though one of them may be a boundary edge).

    Then the two tetrahedra $\Delta$ and $\Delta'$ may be merged into a single
    tetrahedron, with the region between edges $g$ and $h$ and vertex
    $A$ flattened to a single face.  This operation is called a
    {\em 2-1 edge move} and is illustrated in the right hand diagram
    of Figure~\ref{fig-eltmove21e}.  This move preserves the underlying
    3-manifold and reduces the number of tetrahedra in its triangulation
    by one.
\end{lemma}

\begin{proof}
    We employ a strategy similar to that used in the proof
    of Lemma~\ref{l-eltmove20e}.  Consider the disc bounded by edges $g$
    and $h$ that slices through both tetrahedra $\Delta$ and $\Delta'$.
    Since edges $g$ and $h$ are distinct and are not both boundary we can
    crush this disc to a single edge without changing the underlying
    3-manifold.  This reduces the region between
    edges $g$ and $h$ and vertex $A$ to a triangular pillow as
    illustrated in the central diagram of Figure~\ref{fig-eltmove21eexpl}.
    The pillow is retriangulated using two tetrahedra and
    the region between edges $g$ and $h$ and vertex $B$ is
    retriangulated using a single tetrahedron with a new internal edge
    of degree one.

    \begin{figure}[htb]
    \psfrag{A}{{\small $A$}} \psfrag{B}{{\small $B$}}
    \psfrag{C}{{\small $C$}} \psfrag{D}{{\small $D$}}
    \psfrag{E}{{\small $E$}}
    \psfrag{e}{{\small $e$}}
    \psfrag{g}{{\small $g$}} \psfrag{h}{{\small $h$}}
    \psfrag{T}{{\small $\Delta$}} \psfrag{S}{{\small $\Delta'$}}
    \centerline{\includegraphics[scale=0.7]{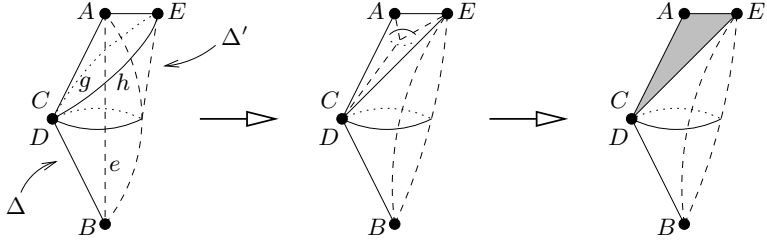}}
    \caption{The intermediate stage of a 2-1 edge move}
    \label{fig-eltmove21eexpl}
    \end{figure}

    Focusing our attention upon the triangular pillow and its interior
    vertex, we can use the constraints upon edges $g$ and $h$ to
    establish the conditions of Lemma~\ref{l-eltmove20v}.
    This allows us to
    perform a 2-0 vertex move upon the pillow, flattening it
    to a single face as illustrated in the right hand diagram of
    Figure~\ref{fig-eltmove21eexpl}.
    This completes the 2-1 edge move with no changes in the
    underlying 3-manifold and an overall reduction of one tetrahedron.
\end{proof}

\subsubsection{Normal Surfaces}

The theory of normal surfaces offers a powerful tool in the development
of algorithms in 3-manifold topology.  Since our use of normal surfaces
does not extend beyond this small section, we refer the reader to
Hemion \cite{hemion} for a detailed overview of normal surface theory
in an algorithmic context.

When constructing a census of 3-manifold triangulations,
an enumeration of the vertex embedded normal surfaces of a
triangulation can help establish whether or not that
triangulation is \ppirr.  In particular, we can call upon
the following results.

\begin{lemma} \label{l-hasnormalplane}
    Let $T$ be a non-orientable 3-manifold triangulation.
    If a projective plane appears amongst the vertex embedded
    normal surfaces of $T$ then $T$ is not {\ppirr}.
\end{lemma}

\begin{proof}
    If the projective plane described above is two-sided, it follows
    immediately from the definition of {\ppirrty}
    that $T$ cannot be {\ppirr}.
    If the projective plane is one-sided on the other hand,
    it follows that a regular neighbourhood of this projective
    plane forms an $\rps$ connected sum component of the underlying
    3-manifold.  Since $\rps$ is orientable but $T$ is
    not, there must be more than one connected
    sum component and so $T$ cannot be irreducible (and therefore $T$
    cannot be \ppirr).
\end{proof}

\begin{lemma} \label{l-hasnoplaneorsphere}
    Let $T$ be a 3-manifold triangulation.
    If the vertex embedded normal surfaces of $T$ contain no projective
    planes and no 2-spheres aside from the trivial
    vertex linking 2-spheres, then $T$ is {\ppirr}.
\end{lemma}

\begin{proof}
    Assume that $T$ is not {\ppirr}.  Then either $T$ contains an
    embedded essential 2-sphere (an embedded 2-sphere that does not
    bound a ball) or $T$ contains an embedded two-sided projective plane.

    It is proven by Kneser \cite{kneser-normal} and
    Schubert \cite{schubert-normal} that any triangulation containing an
    embedded essential 2-sphere contains an embedded essential normal 2-sphere.
    An analogous argument shows that if a triangulation contains
    an embedded projective plane then it contains an embedded
    normal projective plane.  Therefore if $T$ is not {\ppirr}
    then it must contain an embedded normal surface whose
    Euler characteristic is strictly positive.

    This normal surface $S$ can be expressed as
    $S = \lambda_1 V_1+\ldots+\lambda_k V_k$ where each $V_i$ is
    a vertex embedded normal surface and each $\lambda_i > 0$.  Since
    Euler characteristic is a linear
    function of normal surfaces, it follows that
    $\chi(S) = \lambda_1 \chi(V_1)+\ldots+\lambda_k \chi(V_k)$ and in
    particular that $\chi(V_i) > 0$ for some $i$.  Thus some
    vertex embedded normal surface $V_i$ is either a 2-sphere or a
    projective plane.

    Finally we observe that $V_i$ cannot be a vertex linking 2-sphere
    since a vertex link alone is not essential and
    the sum of a vertex link with any other normal surface is disconnected.
\end{proof}

Note that the only case in which we cannot conclude whether or not a
triangulation is {\ppirr} is the case in which its vertex embedded
normal surfaces include a non-vertex-linking 2-sphere but no
projective planes.

\subsubsection{Special Subcomplexes}

There are some particular subcomplexes whose presence within a
triangulation indicates that the triangulation cannot be {\ppirr}.
Two such subcomplexes, the pillow 2-sphere and the snapped 2-sphere, are
described in detail below.

\begin{defn}[Pillow 2-Sphere]
    A {\em pillow 2-sphere} is a 2-sphere formed from two faces of a
    triangulation.  These two faces must be joined along all three
    edges of each face.  Furthermore, these three edges must be
    distinct.  The formation of a pillow 2-sphere is illustrated in
    Figure~\ref{fig-pillow2sphere}.

    \begin{figure}[htb]
    \centerline{\includegraphics[scale=0.7]{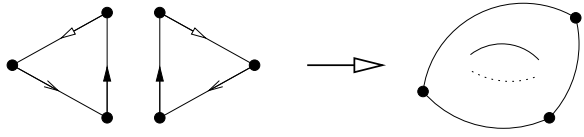}}
    \caption{Forming a pillow 2-sphere}
    \label{fig-pillow2sphere}
    \end{figure}
\end{defn}

\begin{lemma}
    If a 3-manifold triangulation contains a pillow 2-sphere then this
    triangulation cannot be both minimal and {\ppirr}.
\end{lemma}

\begin{proof}
    This result is proven in \cite{burton-facegraphs}.  The proof
    essentially involves converting the pillow 2-sphere to an embedded
    2-sphere, which must bound a ball if the triangulation is to be
    {\ppirr}.  The ball is removed and replaced with a simpler
    structure and the triangulation is then simplified.  The result is a
    triangulation of the same 3-manifold using fewer than the original
    number of tetrahedra, showing the original triangulation to be
    non-minimal.
\end{proof}

\begin{defn}[Snapped 2-Sphere] \label{d-snapped2sphere}
    A {\em snapped 2-sphere} is a 2-sphere formed using two tetrahedra
    $\Delta_1$ and $\Delta_2$ as follows.  Each tetrahedron $\Delta_i$
    is folded upon itself to form an edge $e_i$ of degree one, as
    illustrated in Figure~\ref{fig-snapped2sphere}.  Let
    $f_i$ denote the edge opposite $e_i$ in tetrahedron $\Delta_i$, so
    that $f_i$ bounds a disc $\delta_i$ slicing through the midpoint of edge
    $e_i$.  Discs $\delta_1$ and $\delta_2$ are shaded in the diagram.

    \begin{figure}[htb]
    \psfrag{T1}{{\small $\Delta_1$}} \psfrag{T2}{{\small $\Delta_2$}}
    \psfrag{D1}{{\small $\delta_1$}} \psfrag{D2}{{\small $\delta_2$}}
    \psfrag{e1}{{\small $e_1$}} \psfrag{e2}{{\small $e_2$}}
    \psfrag{f1}{{\small $f_1$}} \psfrag{f2}{{\small $f_2$}}
    \centerline{\includegraphics[scale=0.7]{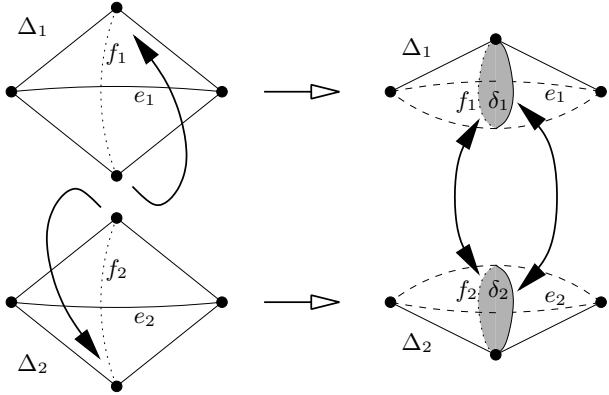}}
    \caption{Forming a snapped 2-sphere}
    \label{fig-snapped2sphere}
    \end{figure}

    Edges $f_1$ and $f_2$ are then identified in the triangulation.
    As a result discs $\delta_1$ and $\delta_2$ are joined at their boundaries
    to form a 2-sphere within the triangulation, referred to
    as a snapped 2-sphere.
\end{defn}

\begin{lemma}
    If a 3-manifold triangulation contains a snapped 2-sphere then this
    triangulation cannot be both minimal and {\ppirr}.
\end{lemma}

\begin{proof}
    Let $T$ be a triangulation of the 3-manifold $M$ and assume that $T$
    is \ppirr.  Furthermore, assume that $T$ contains a snapped 2-sphere
    as described in Definition~\ref{d-snapped2sphere}.

    Note that a snapped 2-sphere is an embedded 2-sphere, since it
    contains only one vertex and one edge and therefore has no possible
    points of self-intersection.  As with all 2-spheres, it is also
    two-sided.
    We can thus slice triangulation $T$ along this 2-sphere to
    obtain two spherical boundary components, each formed from a copy
    of the two discs $\delta_1$ and $\delta_2$.

    We can now effectively cap each of these boundary
    components with a ball by identifying its two boundary discs together,
    producing the new triangulation $T'$ as illustrated in
    Figure~\ref{fig-snapped2spheremixed}.  Note that $T$ and $T'$
    contain precisely the same number of tetrahedra.

    \begin{figure}[htb]
    \psfrag{T1}{{\small $\Delta_1$}} \psfrag{T2}{{\small $\Delta_2$}}
    \psfrag{D1}{{\small $\delta_1$}} \psfrag{D2}{{\small $\delta_2$}}
    \psfrag{e1}{{\small $e_1$}} \psfrag{e2}{{\small $e_2$}}
    \psfrag{f1}{{\small $f_1$}} \psfrag{f2}{{\small $f_2$}}
    \psfrag{Tri1}{{\small Triangulation $T$}}
    \psfrag{Tri2}{{\small Triangulation $T'$}}
    \centerline{\includegraphics[scale=0.7]{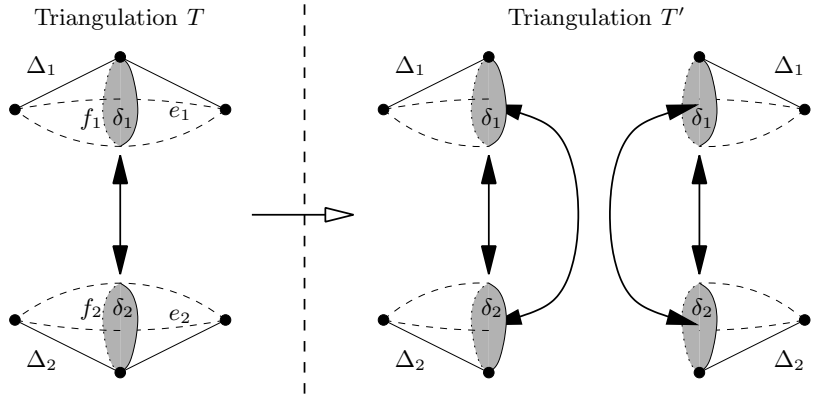}}
    \caption{Slicing along a snapped 2-sphere and capping the boundary spheres}
    \label{fig-snapped2spheremixed}
    \end{figure}

    Since $T$ is {\ppirr}, the snapped 2-sphere must be separating in
    $T$ and so triangulation $T'$ consists of two disconnected
    components.  Moreover, these two components represent a connected
    sum decomposition of $M$.
    Again by {\ppirrty}, it follows that one of these
    components is a 3-sphere and the other is a new triangulation of
    the original 3-manifold $M$.  Since this new triangulation of $M$
    contains strictly fewer tetrahedra than $T$, it follows
    that the original triangulation $T$ cannot be minimal.
\end{proof}

\subsubsection{Invariant Analysis} \label{s-algmanalysis-inv}

Once a triangulation is known to be {\ppirr}, invariant analysis can
be used to prove its minimality.  This technique requires the census
to be constructed according to
increasing numbers of tetrahedra, i.e., all
1-tetrahedron triangulations should be generated and analysed, then
all 2-tetrahedron triangulations, all 3-tetrahedron triangulations
and so on.

The key observation is that if some {\ppirr} triangulation $T$ is
non-minimal, then a triangulation of the same 3-manifold must appear
in an earlier section of the census formed from
fewer tetrahedra.  Thus, if some
collection of 3-manifold invariants can be found that together
distinguish the underlying 3-manifold of $T$ from any of the 3-manifolds
constructed in earlier sections of the census, it follows that $T$ is a
minimal {\ppirr} triangulation and is thus eligible to appear in our
final list of census results.

The invariants that were used for this particular census include homology
groups, fundamental group and the quantum invariants of
Turaev and Viro \cite{turaev-viro}.  The Turaev-Viro invariants,
used also by Matveev \cite{matveev9}, have proven exceptionally
useful in practice for distinguishing 3-manifolds in both orientable
and non-orientable settings.

\subsubsection{Results} \label{s-algmanalysis-results}

For the $\leq 7$-tetrahedron non-orientable census described in this paper,
the techniques discussed in Sections~\ref{s-algmanalysis-eltmoves}
to~\ref{s-algmanalysis-inv} prove sufficient to determine precisely
which triangulations are minimal and {\ppirr}.  Furthermore, a
combination of elementary moves and invariant analysis allows
these minimal {\ppirr} triangulations to be grouped
into equivalence classes according to which triangulations have
identical underlying 3-manifolds.  In this way
we obtain a final census of 41 triangulations
representing eight different 3-manifolds.

We proceed now to Section~\ref{s-common} which paves the way for a
combinatorial analysis of these 41 census triangulations.

\section{Common Structures} \label{s-common}

In order to make the census triangulations easier to both visualise and
analyse, we decompose these triangulations into a variety of building
blocks.  Ideally such building blocks should be large enough that they
significantly simplify the representation and analysis of the
triangulations containing them, yet small enough that they can be
frequently reused throughout the census.

This idea of describing triangulations using medium-sized building
blocks has been used previously for the orientable case.  Matveev
describes a few orientable building blocks \cite{matveev6} and Martelli
and Petronio describe a more numerous set of smaller
orientable building blocks called {\em bricks} \cite{italian9}.

An examination of the non-orientable triangulations of this census
shows a remarkable consistency of combinatorial structure.
We therefore need only two types of building block:
the thin $I$-bundle and the layered solid torus.

\subsection{Thin $I$-Bundles} \label{s-thin}

A thin $I$-bundle is essentially a triangulation of
an $I$-bundle that has a thickness of
only one tetrahedron between its two parallel boundaries.  Thin
$I$-bundles play an important role in the construction of minimal
non-orientable triangulations and appear within all of the triangulations
described in this census.

A thin $I$-bundle over a surface $S$ is built upon a decomposition of
$S$ into triangles and quadrilaterals.  Each triangle or quadrilateral
of $S$ corresponds to a single tetrahedron of the overall $I$-bundle.

We thus begin by discussing triangle and quadrilateral decompositions of
surfaces and the properties that we require of such decompositions.

\begin{defn}[Well-Balanced Decomposition] \label{d-wellbalanced}
    Let $S$ be some closed surface.  A {\em well-bal\-anced decomposition}
    of $S$ is a decomposition of $S$ into triangles and quadrilaterals
    satisfying the following properties.
    \begin{enumerate}
        \item \label{en-wb-parity}
        Every vertex of the decomposition meets an even number of
        quadrilateral corners.
        \item \label{en-wb-tridiscs}
        If the quadrilaterals are removed then the surface breaks
        into a disconnected collection of triangulated discs.
        \item \label{en-wb-quadcycles}
        There are no cycles of quadrilaterals.  That is, any path
        formed by walking through a series of quadrilaterals, always
        entering and exiting by opposite sides, must eventually run into
        a triangle.
    \end{enumerate}
\end{defn}

\begin{figure}[htb]
\centerline{\includegraphics[scale=0.7]{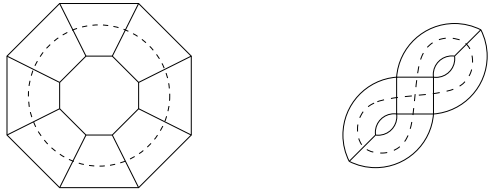}}
\caption{Arrangements of quadrilaterals that include cycles}
\label{fig-quadcyclesbad}
\end{figure}

\begin{figure}[htb]
\centerline{\includegraphics[scale=0.7]{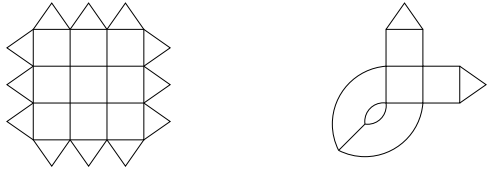}}
\caption{Arrangements of quadrilaterals that do not include cycles}
\label{fig-quadcyclesgood}
\end{figure}

To clarify condition~\ref{en-wb-quadcycles},
Figure~\ref{fig-quadcyclesbad} illustrates
some arrangements of quadrilaterals that contain cycles as described
above.  In each diagram the offending cycle is
marked by a dotted line.  Figure~\ref{fig-quadcyclesgood} on the other
hand illustrates arrangements of quadrilaterals that do not contain
cycles and so are perfectly acceptable within a well-balanced
decomposition.

\begin{example}
    Figure~\ref{fig-wellbalancedgood} illustrates a handful of well-balanced
    decompositions of the torus.  Each of these decompositions can be
    seen to satisfy all of the conditions of Definition~\ref{d-wellbalanced}.

    \begin{figure}[htb]
    \centerline{\includegraphics[scale=0.7]{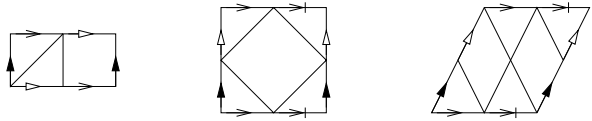}}
    \caption{Well-balanced decompositions of the torus}
    \label{fig-wellbalancedgood}
    \end{figure}

    Figure~\ref{fig-wellbalancedbad} however illustrates some
    decompositions of the torus that are not well-balanced.  The first
    diagram illustrates a decomposition that breaks
    condition~\ref{en-wb-parity} of Definition~\ref{d-wellbalanced}; one
    of the offending vertices is marked with a black circle.  The second
    diagram shows how this same decomposition breaks
    condition~\ref{en-wb-quadcycles}; two cycles of quadrilaterals are
    marked with dotted lines.  The final decomposition breaks
    condition~\ref{en-wb-tridiscs} (amongst others); the shaded
    triangles together form an annulus, not a disc.

    \begin{figure}[htb]
    \centerline{\includegraphics[scale=0.7]{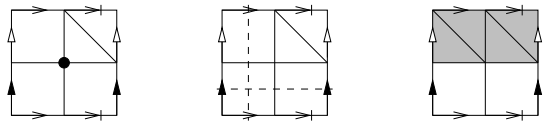}}
    \caption{Decompositions of the torus that are not well-balanced}
    \label{fig-wellbalancedbad}
    \end{figure}
\end{example}

Once we have obtained a well-balanced decomposition of a surface, we
can flesh out its triangles and quadrilaterals
into a full 3-manifold triangulation as follows.

\begin{defn}[Enclosing Triangulation]
    Let $D$ be a well-balanced decomposition of a closed surface.
    The {\em enclosing triangulation} of $D$ is the unique triangulation
    formed as follows.

    Each triangle or quadrilateral of $D$ is enclosed within its own
    tetrahedron as illustrated in Figure~\ref{fig-enclosediscs}.  The
    faces of these tetrahedra are then identified in the unique manner
    that causes these discs to be connected according to the decomposition
    $D$.  Any tetrahedron faces that do not meet the discs of $D$
    (i.e., the faces parallel to triangular discs) are left to
    become boundary faces of the triangulation.

    \begin{figure}[htb]
    \centerline{\includegraphics[scale=0.7]{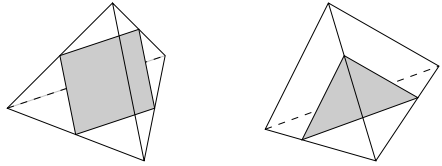}}
    \caption{Enclosing triangles and quadrilaterals within tetrahedra}
    \label{fig-enclosediscs}
    \end{figure}

    As an example, Figure~\ref{fig-splitreconstruct} illustrates two
    quadrilaterals $q$ and $q'$ placed within
    within tetrahedra $\Delta$ and $\Delta'$ respectively.
    Suppose edges $e$ and $e'$ of these quadrilaterals are
    identified within the surface decomposition $D$.
    Then the tetrahedron faces $f$ and $f'$ must be identified as
    illustrated so that edges $e$ and $e'$ can be connected correctly.

    \begin{figure}[htb]
    \psfrag{D1}{{\small $\Delta$}} \psfrag{D2}{{\small $\Delta'$}}
    \psfrag{q1}{{\small $q$}} \psfrag{q2}{{\small $q'$}}
    \psfrag{e1}{{\small $e$}} \psfrag{e2}{{\small $e'$}}
    \psfrag{f1}{{\small $f$}} \psfrag{f2}{{\small $f'$}}
    \centerline{\includegraphics[scale=0.7]{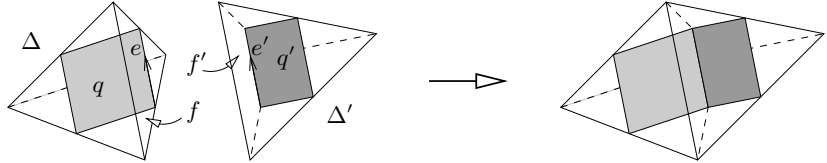}}
    \caption{Constructing the enclosing triangulation}
    \label{fig-splitreconstruct}
    \end{figure}
\end{defn}

Since each disc of the decomposition $D$ runs through the centre
of a tetrahedron,
we see that the enclosing triangulation simply thickens the
2-manifold decomposition $D$ into a 3-dimensional space.  The
original surface with its triangles and quadrilaterals in turn becomes
a normal surface running through the centre of the enclosing triangulation.

\begin{lemma} \label{l-thinibundle}
    Let $D$ be a well-balanced decomposition of the closed surface $S$.
    Then the enclosing triangulation of $D$ is a 3-manifold triangulation
    representing an $I$-bundle (possibly twisted) over $S$.
\end{lemma}

\begin{proof}
    Let $T$ be the enclosing triangulation.  Our first step
    is to prove that $T$ is in fact a triangulation of a 3-manifold.

    It can be shown that each edge not meeting the central decomposition
    $D$ is in fact a boundary edge of the triangulation.  This is
    clearly true of edges belonging to tetrahedra that enclose triangles
    of $D$, since faces parallel to triangles of $D$
    become boundary faces of the overall triangulation.
    This is also true of edges belonging to tetrahedra
    that enclose quadrilaterals of $D$, since
    condition~\ref{en-wb-quadcycles} of Definition~\ref{d-wellbalanced}
    (no cycles of quadrilaterals) ensures that each such edge is
    identified with an edge of one of the boundary faces previously described.

    Each internal edge of $T$ therefore cuts through a vertex of
    the decomposition $D$.  Condition~\ref{en-wb-parity} of
    Definition~\ref{d-wellbalanced} (vertices meet an even number of
    quadrilateral corners) ensures that no internal edge of
    $T$ is identified with itself in reverse.

    We turn now to the vertices of $T$.  By observing how
    the tetrahedra are hooked together we see that the link of each
    vertex $v$ is of the form illustrated in
    Figure~\ref{fig-thinvertexlink}.

    \begin{figure}[htb]
    \centerline{\includegraphics[scale=0.7]{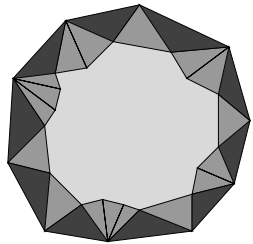}}
    \caption{A vertex link in an enclosing triangulation}
    \label{fig-thinvertexlink}
    \end{figure}

    The dark shaded triangles around the
    outside correspond to tetrahedra enclosing triangles of
    $D$ in which $v$ belongs to a boundary face.  This is illustrated in the
    left hand diagram of Figure~\ref{fig-thinvertexlinkpieces}.
    These dark shaded triangles provide the boundary edges of the
    vertex link.

    \begin{figure}[htb]
    \psfrag{D}{{\small $D$}}
    \psfrag{v}{{\small $v$}}
    \psfrag{Link}{{\small Link}}
    \centerline{\includegraphics[scale=0.7]{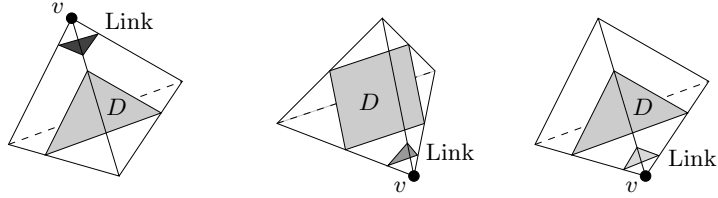}}
    \caption{Pieces of a vertex link}
    \label{fig-thinvertexlinkpieces}
    \end{figure}

    The medium shaded triangles just inside this boundary correspond to
    tetrahedra enclosing quadrilaterals of $D$, as illustrated in
    the central diagram of Figure~\ref{fig-thinvertexlinkpieces}.  Each
    of these medium shaded triangles meets the boundary at a single
    point, since we observed earlier that each edge of $T$
    running parallel to a quadrilateral of $D$ is in fact a boundary edge.

    The remaining pieces of the vertex link, corresponding to the light
    shaded region in the interior of Figure~\ref{fig-thinvertexlink},
    are provided by
    tetrahedra enclosing triangles of $D$ in which $v$ lies opposite
    the boundary face.  In these cases the pieces of the vertex link
    run parallel to the triangles of $D$
    as illustrated in the right hand diagram of
    Figure~\ref{fig-thinvertexlinkpieces}.  From
    condition~\ref{en-wb-tridiscs} of Definition~\ref{d-wellbalanced}
    (triangulated regions form discs in $D$) we see that these interior
    pieces combine to form a topological disc.  Thus $v$ is a
    boundary vertex (i.e., the darker band forming the boundary of
    Figure~\ref{fig-thinvertexlink}
    actually exists), and more importantly the entire link of
    vertex $v$ is a topological disc.

    Thus our enclosing triangulation $T$ is indeed a triangulation of a
    3-ma\-ni\-fold.  From the fact that the well-balanced
    decomposition of the surface $S$
    runs through the centre of each tetrahedron, as well as the earlier
    observation that each vertex, edge and face not meeting this surface
    decomposition in fact forms part of the triangulation boundary, it is
    straightforward to see that the enclosing triangulation represents
    an $I$-bundle (possibly twisted) over $S$.
\end{proof}

\begin{defn}[Thin $I$-Bundle] \label{d-thinibundle}
    A {\em thin $I$-bundle} over a closed surface $S$ is the
    enclosing triangulation of a well-balanced decomposition of $S$.
    This well-balanced decomposition is referred to as the {\em central
    surface decomposition} of the thin $I$-bundle.
\end{defn}

We see then that Lemma~\ref{l-thinibundle} simply states that a thin
$I$-bundle over a surface $S$ is what it claims to be, i.e., an actual
triangulation of an $I$-bundle over $S$.

Before closing this section we present a method of marking
a well-balanced decomposition that allows us to establish precisely how
the corresponding $I$-bundle is twisted, if at all.

\begin{defn}[Marked Decomposition] \label{d-markings}
    A well-balanced decomposition can be {\em marked} to illustrate how
    it is embedded within its enclosing triangulation.  Markings consist
    of solid lines and dotted lines, representing features above and
    below the central surface respectively.

    The different types of markings are illustrated in
    Figure~\ref{fig-markings}.  Each quadrilateral lies between two
    perpendicular boundary edges of the triangulation, one above and one
    below.  These boundary edges are represented by a solid line and a
    dotted line as illustrated in the left hand diagram of
    Figure~\ref{fig-markings}.  Each triangle lies between a boundary
    vertex and a boundary face.  If the boundary vertex lies above the
    triangle (and the boundary face below) then the triangle is marked
    with three solid lines as illustrated in the central diagram of
    Figure~\ref{fig-markings}.  If the boundary vertex lies below the
    triangle (and the boundary face above) then the triangle is marked
    with three dotted lines as illustrated in the right hand diagram.

    \begin{figure}[htb]
    \centerline{\includegraphics[scale=0.7]{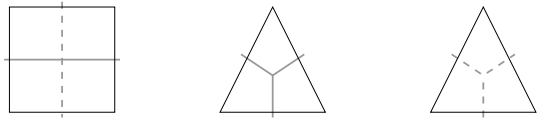}}
    \caption{Marking quadrilaterals and triangles}
    \label{fig-markings}
    \end{figure}

    Since the edges and vertices of the $I$-bundle boundary all connect
    together, it follows that the markings must similarly connect
    together, with solid lines matched with solid lines and dotted lines
    matched with dotted lines.  Following these markings across edge
    identifications therefore allows us to see if and where a region above the
    central surface moves through a twist to become a region below.
\end{defn}

\begin{example}
    Figure~\ref{fig-wellbalancedmarked} illustrates three well-balanced
    decompositions of the torus with markings.  In the first diagram we
    see that the $I$-bundle is twisted across the upper and lower
    edge identifications, since solid lines change to dotted lines across
    these identifications and vice versa.  There is no twist however
    across the left and right edge identifications since the
    solid line does not change.

    \begin{figure}[htb]
    \centerline{\includegraphics[scale=0.7]{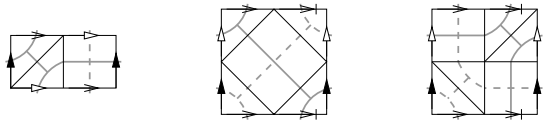}}
    \caption{Markings on well-balanced torus decompositions}
    \label{fig-wellbalancedmarked}
    \end{figure}

    In the second diagram we see that the $I$-bundle
    is twisted across all of the outer edge identifications, with solid
    and dotted lines being exchanged in every case.
    In the third diagram we see that there are no twists at all, and that
    the $I$-bundle is in fact simply the product $\torus \times I$.
\end{example}

\subsection{Layered Solid Tori} \label{s-lst}

A key structure that appears frequently within both orientable and
non-ori\-ent\-able minimal triangulations is the layered solid torus.
Layered solid tori have been discussed in a variety of informal contexts
by Jaco and Rubinstein.  They appear in \cite{0-efficiency}
and are treated thoroughly in \cite{layeredlensspaces}.
Analogous constructs involving special spines of
3-manifolds can be found in \cite{italianfamilies} and \cite{matveev6}.
The preliminary definitions presented here follow those given in
\cite{burton-facegraphs}.

In order to describe the construction of a layered solid torus we
introduce
the process of layering.
Layering is a transformation that, when applied to a triangulation with
boundary, does not change the underlying 3-manifold but does change the
curves formed by the boundary edges of the triangulation.

\begin{defn}[Layering] \label{d-layering}
    Let $T$ be a triangulation with boundary and let $e$ be one of its
    boundary edges.  To {\em layer a tetrahedron on edge $e$}, or just
    to {\em layer on edge $e$}, is to take
    a new tetrahedron $\Delta$, choose two of its faces and identify them
    with the two boundary faces on either side of
    $e$ without twists.  This procedure is
    illustrated in Figure~\ref{fig-layering}.

    \begin{figure}[htb]
    \psfrag{e}{{\small $e$}}
    \psfrag{f}{{\small $f$}}
    \psfrag{T}{{\small $\Delta$}}
    \psfrag{M}{{\small $T$}}
    \centerline{\includegraphics[width=5.5cm]{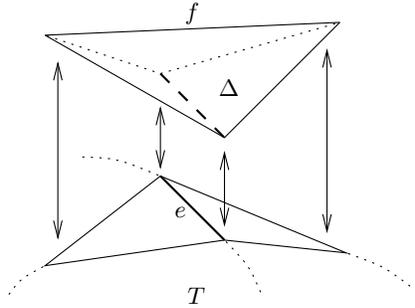}}
    \caption{Layering a tetrahedron on a boundary edge}
    \label{fig-layering}
    \end{figure}
\end{defn}

Note that layering on a boundary edge does not change the underlying
3-manifold; the only effect is to thicken the boundary around edge
$e$.  Moreover, once a layering has been performed, edge $e$ is no
longer a boundary edge but instead edge $f$
(which in general represents a different curve on
the boundary of the 3-manifold) has been added as a new boundary edge.

\begin{defn}[Layered Solid Torus] \label{d-lst}
    A {\em standard layered solid torus} is a triangulation of a solid torus
    formed as follows.  We begin with the one-triangle
    {\mobius} band illustrated in the left hand diagram of
    Figure~\ref{fig-layermobius}, where the two edges marked $e$ are
    identified according to the arrows and where $g$ is a boundary edge.
    This {\mobius} band can be embedded in $\R^3$ as
    illustrated in the right hand diagram of Figure~\ref{fig-layermobius}.

    \begin{figure}[htb]
    \psfrag{e}{{\small $e$}} \psfrag{f}{{\small $f$}} \psfrag{g}{{\small $g$}}
    \psfrag{F}{{\small $F$}} \psfrag{F'}{{\small $F^\prime$}}
    \centerline{\includegraphics[scale=0.7]{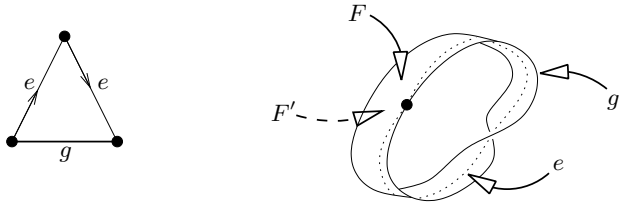}}
    \caption{A one-triangle \mobius\ band}
    \label{fig-layermobius}
    \end{figure}

    In this embedding our single triangular face has two sides, marked
    $F$ and $F'$ in the diagram.  We make an initial layering upon edge
    $e$ as illustrated in Figure~\ref{fig-layermobiusnewtet}, so that
    faces {\em ABC} and {\em BCD} of the new tetrahedron are joined to
    sides $F$ and $F'$ respectively of the original triangular face.
    Although the initial {\mobius} band is not actually a 3-manifold
    triangulation, the layering procedure remains the same as described
    in Definition~\ref{d-layering}.

    \begin{figure}[htb]
    \psfrag{A}{{\small $A$}} \psfrag{B}{{\small $B$}}
    \psfrag{C}{{\small $C$}} \psfrag{D}{{\small $D$}}
    \psfrag{e}{{\small Layer upon edge $e$}}
    \psfrag{F}{{\small $F$}} \psfrag{F'}{{\small $F^\prime$}}
    \centerline{\includegraphics[scale=0.7]{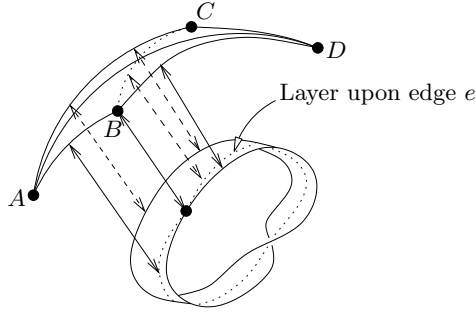}}
    \caption{Layering on a {\mobius} band to form a solid torus}
    \label{fig-layermobiusnewtet}
    \end{figure}

    Since $F$ and $F'$ are in fact opposite sides of the same triangular
    face, we see that faces {\em ABC} and {\em BCD} become identified as
    illustrated in Figure~\ref{fig-layermobiusstandard}.  The result is
    the well-known one-tetrahedron triangulation of the solid torus.
    The identified faces {\em ABC} and {\em BCD} are shaded in the diagram.

    \begin{figure}[htb]
    \psfrag{A}{{\small $A$}} \psfrag{B}{{\small $B$}}
    \psfrag{C}{{\small $C$}} \psfrag{D}{{\small $D$}}
    \psfrag{e}{{\small $e$}} \psfrag{g}{{\small $g$}}
    \psfrag{F}{{\small $F$}} \psfrag{F'}{{\small $F^\prime$}}
    \centerline{\includegraphics[scale=0.7]{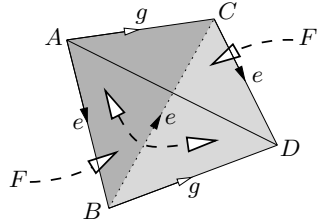}}
    \caption{A one-tetrahedron standard layered solid torus}
    \label{fig-layermobiusstandard}
    \end{figure}

    Having obtained a 3-manifold triangulation of the solid
    torus, we finish the construction by performing some number of
    additional layerings upon boundary edges, one at a time.  We may
    layer as many times we like or we may make no additional layerings
    at all.  There are thus infinitely many different standard layered
    solid tori that can be constructed.

    It is useful to consider the {\mobius} band on its own
    to be a {\em degenerate layered solid torus} containing zero tetrahedra.
    A {\em non-standard layered solid torus} can also be formed by
    making the initial layering upon edge $g$ of the {\mobius} band
    instead of edge $e$, although such structures are not considered here.
\end{defn}

We can observe that each standard layered solid torus has two boundary faces
and represents the same underlying 3-manifold, i.e., the solid torus.
What distinguishes the different layered solid tori
is the different patterns of curves that their boundary edges
make upon the boundary torus.

\begin{defn}[Three-Parameter Torus Curves] \label{d-toruscurve}
    Let $T$ be a torus formed from two triangles as illustrated in
    Figure~\ref{fig-layertorusplusminus}.  Label one of these triangles
    $+$ and the other $-$.  Select some ordering of the three edges
    and label these edges $e_1$, $e_2$ and $e_3$ accordingly.

    \begin{figure}[htb]
    \psfrag{+}{{\small $+$}} \psfrag{-}{{\small $-$}}
    \psfrag{e1}{{\small $e_1$}} \psfrag{e2}{{\small $e_2$}}
    \psfrag{e3}{{\small $e_3$}}
    \centerline{\includegraphics{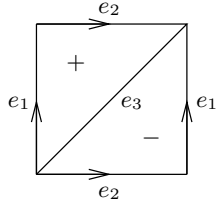}}
    \caption{A two-triangle torus}
    \label{fig-layertorusplusminus}
    \end{figure}

    Consider some oriented closed curve on this torus
    as illustrated in Figure~\ref{fig-layer235curve}.  Using this curve
    we can assign a number to each edge $e_i$, this being
    the number of
    times the curve crosses edge $e_i$ from $+$ to $-$ minus the number of
    times it crosses edge $e_i$ from $-$ to $+$.

    \begin{figure}[htb]
    \psfrag{+}{{\small $+$}} \psfrag{-}{{\small $-$}}
    \psfrag{e1}{{\small $e_1$}} \psfrag{e2}{{\small $e_2$}}
    \psfrag{e3}{{\small $e_3$}}
    \centerline{\includegraphics{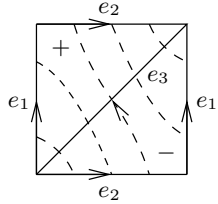}}
    \caption{A $(2,3,-5)$ curve on a torus}
    \label{fig-layer235curve}
    \end{figure}

    If the numbers assigned to edges $e_1$, $e_2$ and $e_3$ are
    $p$, $q$ and $r$ respectively, we refer to our oriented curve as a
    {\em $(p,q,r)$ curve}.  Thus, for instance, the curve illustrated in
    Figure~\ref{fig-layer235curve} is a $(2,3,-5)$ curve.
\end{defn}

It is trivial to show that any $(p,q,r)$ curve satisfies $p+q+r=0$.  It
is also straightforward to show that that if the curve is embedded then
either $p$, $q$ and $r$ are pairwise coprime or $(p,q,r)=(0,0,0)$.
We can use three-parameter torus curves to categorise layered solid tori as
follows.

\begin{defn}[Layered Solid Torus Parameters] \label{d-lstparams}
    Let $L$ be a layered solid torus.  Upon the two faces that form the
    boundary torus, draw the boundary of a meridinal disc of the
    underlying solid torus.  Assign to this meridinal curve some arbitrary
    orientation and arbitrarily label the two boundary faces $+$ and $-$.

    The meridinal curve then forms some $(p,q,r)$ curve on the boundary
    torus.  $p$, $q$ and $r$ are said to be the {\em parameters} of
    the layered solid torus $L$, and $L$ is said to be a
    {\em $(p,q,r)$ layered solid torus}, denoted $\lst(p,q,r)$.
\end{defn}

Note that a $(p,q,r)$ layered solid torus is also a $(-p,-q,-r)$ layered
solid torus.  Note furthermore that the layered solid torus parameters are not
ordered, so for instance a $(p,q,r)$ layered solid torus could also be
written with parameters $(p,r,q)$ or $(q,r,p)$.

\begin{example} \label{ex-toruslst123}
    A $(1,2,-3)$ layered solid torus formed from one tetrahedron is
    illustrated in Figure~\ref{fig-layerlst123}.  This is in fact the
    same layered solid torus illustrated in
    Figure~\ref{fig-layermobiusstandard}, formed by a single
    layering upon edge $e$ of the {\mobius} band of
    Figure~\ref{fig-layermobius}.

    \begin{figure}[htb]
    \psfrag{P}{{\small $P$}} \psfrag{Q}{{\small $Q$}}
    \psfrag{R}{{\small $R$}} \psfrag{S}{{\small $S$}}
    \centerline{\includegraphics[scale=0.7]{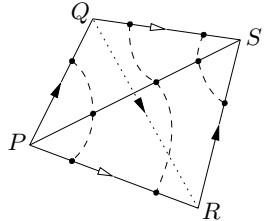}}
    \caption{A $(1,2,-3)$ layered solid torus}
    \label{fig-layerlst123}
    \end{figure}

    The back two faces of the tetrahedron are identified;
    specifically face $PQR$ is identified with face $QRS$.
    The meridinal curve is illustrated by the dotted line drawn upon
    the boundary faces, and its intersections with
    the edges of the boundary torus are marked.
\end{example}

\begin{example}
    The degenerate layered solid torus with zero tetrahedra, i.e., the
    {\mobius} band, is an example of
    a $(1,1,-2)$ layered solid torus.  The meridinal
    curve is marked in Figure~\ref{fig-layerlst112} with a dotted line.
    The three edges of the (degenerate) boundary
    torus are $g$ and the front and rear sides of $e$ (recalling that
    the {\mobius} band is embedded in $\R^3$).
    Once more the intersections of the meridinal curve with these
    boundary edges are marked with black circles in the diagram.

    \begin{figure}[htb]
    \psfrag{e}{{\small $e$}} \psfrag{f}{{\small $f$}} \psfrag{g}{{\small $g$}}
    \centerline{\includegraphics{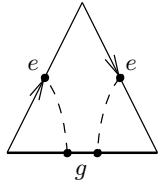}}
    \caption{The degenerate $\lst(1,1,-2)$}
    \label{fig-layerlst112}
    \end{figure}
\end{example}

Recall from Definition~\ref{d-lst} that a layered solid torus is built
by layering upon boundary edges one at a time.  We can trace the
parameters of the layered solid torus as these successive layerings take
place.  The following result is relatively straightforward to prove by
keeping track of intersections with boundary edges.

\begin{lemma} \label{l-toruslstgrowth}
    Layering on the edge with parameter $r$ in a $(p,q,r)$ layered solid
    torus produces a $(p,-q,q-p)$ layered solid torus.
\end{lemma}

Armed with this result we can present a general method for constructing a
layered solid torus with a given set of parameters.

\begin{algorithm}[Layered Solid Torus Construction] \label{a-lstconstruction}
    For any pairwise coprime integers $p$, $q$ and $r$ satisfying
    $p+q+r=0$ and $\max(|p|,|q|,|r|) \geq 3$, the following algorithm
    can be used to construct a standard layered solid torus with
    parameters $(p,q,r)$.  Without loss of generality, let
    $|r| \geq |p|,|q|$.

    Since $p$, $q$ and $r$ are pairwise
    coprime we see that in fact $|r| > |p|,|q|$, and since
    $p+q+r=0$ we also see that $|p|,|q| > |q-p|$.
    If either $|p| \geq 3$ or $|q| \geq 3$ we can therefore use
    this same algorithm to construct a smaller $\lst(p, -q, q-p)$ and then
    layer on the edge with parameter $|q-p|$ to obtain an $\lst(p,q,r)$.
    Since $\max(|p|,|q|,|r|)$
    decreases at each stage we are guaranteed that this recursion will
    eventually terminate.

    If $|p|,|q| < 3$ then the only triples satisfying the given
    constraints are $(p,q,r) = \pm(1,2,-3)$.  In this case we use the
    one-tetrahedron $\lst(1,2,-3)$ illustrated in
    Example~\ref{ex-toruslst123}.
\end{algorithm}

\begin{example}
    Suppose we wish to construct an $\lst(3,7,-10)$ using
    Algorithm~\ref{a-lstconstruction}.
    To form an $\lst(3,7,-10)$ we layer on edge $4$ in an $\lst(3,4,-7)$.
    To form an $\lst(3,4,-7)$ we layer on edge $1$ in an $\lst(1,3,-4)$.
    To form an $\lst(1,3,-4)$ we layer on edge $2$ in an $\lst(1,2,-3)$,
    and an $\lst(1,2,-3)$ is given to us in Example~\ref{ex-toruslst123}.
    The resulting triangulation has four tetrahedra.
\end{example}

Algorithm~\ref{a-lstconstruction} is in fact the optimal method for
constructing a layered solid torus, as shown by the following result of
Jaco and Rubinstein \cite{layeredlensspaces}.

\begin{theorem}
    Let $p$, $q$ and $r$ be pairwise coprime integers for which
    $p+q+r=0$ and $\max(|p|,|q|,|r|) \geq 3$.  Then there is a unique
    $\lst(p,q,r)$ up to isomorphism that uses the least possible number
    of tetrahedra, and this is the $\lst(p,q,r)$ constructed by
    Algorithm~\ref{a-lstconstruction}.
\end{theorem}

Throughout the following sections, any reference to an $\lst(p,q,r)$ is
assumed to refer to this unique minimal $\lst(p,q,r)$.

\section{Families of Closed Triangulations} \label{s-families}

Having developed a set of medium-sized building blocks in
Section~\ref{s-common}, we can now combine these building blocks
to form closed triangulations.  In this section we present a
number of families of closed triangulations, each of which represents an
infinite class of triangulations sharing a common large-scale structure.
It is seen in Section~\ref{s-census} that, with the exception of the three
triangulations described in Section~\ref{s-exceptional},
the four families presented here together encompass all closed non-orientable
minimal {\ppirr} triangulations formed from up to seven tetrahedra.

A categorisation of triangulations into infinite families as described
above is certainly appealing.  Large classes of
triangulations may be studied simultaneously, and algorithms
become available for generating triangulations of infinite
classes of 3-manifolds.

Furthermore, when presented with an arbitrary
triangulation of an unknown 3-manifold, having
a rich collection of such families at our disposal increases the
chance that we can identify this 3-manifold.
Specifically, if we can manipulate the triangulation into a form that is
recognised as a member of one of these families then the underlying
3-manifold can be subsequently established.

In the case of orientable 3-manifolds, several infinite
parameterised families of triangulations are described in the
literature \cite{burton-thesis,italian9,italianfamilies,matveev6}.

\subsection{Notation} \label{s-defsurfacebundles}

We begin by outlining some notation that is used to describe torus and
Klein bottle bundles over the circle.

\begin{itemize}
    \item $\torus \times I / A$ represents a torus
    bundle over the circle, where $A$ is a unimodular
    $2 \times 2$ matrix indicating the
    homeomorphism under which the torus $\torus \times 0$ is identified
    with the torus $\torus \times 1$.  Note that this space is
    orientable or non-orientable according to whether the determinant
    of $A$ is $+1$ or $-1$.

    More specifically, let $\mu_0$ and $\lambda_0$ be closed curves that
    together generate the fundamental group of the first torus and let
    $\mu_1$ and $\lambda_1$ be the curves parallel to these on the second
    torus.  If
    \[ A = \homtwolarge{a}{b}{c}{d}, \]
    then the homeomorphism under which the two
    tori are identified maps curve $\mu_0$ to $\mu_1^a \lambda_1^c$ and
    curve $\lambda_0$ to $\mu_1^b \lambda_1^d$.

    \item $\kb \times I / A$ represents a Klein bottle
    bundle over the circle, where $A$ is again a unimodular
    $2 \times 2$ matrix indicating the
    homeomorphism under which the Klein bottle $\kb \times 0$ is identified
    with the Klein bottle $\kb \times 1$.

    Let $\mu_0$ and $\lambda_0$ be orientation-preserving and
    orientation-reversing closed curves respectively on the first Klein
    bottle that meet transversely in a single
    point.  It is known that every element of the
    fundamental group of this Klein bottle can be represented as
    $\mu^p \lambda^q$ for some unique pair of integers $p$ and $q$.

    Let $\mu_1$ and $\lambda_1$ be the corresponding curves on the
    second Klein bottle.  If
    \[ A = \homtwolarge{a}{b}{c}{d}, \]
    then the homeomorphism under which the two Klein bottles are
    identified maps curve $\mu_0$ to $\mu_1^a \lambda_1^c$ and
    curve $\lambda_0$ to $\mu_1^b \lambda_1^d$.  It is shown in
    \cite{hempel} that every such matrix $A$ must be of the form
    \[ A = \homtwolarge{\pm1}{b}{0}{\pm1} . \]
\end{itemize}

The notation described above is consistent with that used by
Matveev for orientable 3-manifolds \cite{matveev6}.

\subsection{Layered Surface Bundles} \label{s-surfacebundles}

Our first family of triangulations is the family of layered surface
bundles, constructed from untwisted thin $I$-bundles as follows.

\begin{defn}[Layered Surface Bundle] \label{d-layeredsurfacebundle}
    A {\em layered surface bundle} is a triangulation of a closed
    surface bundle over the circle formed using the
    following construction.  We begin with a thin $I$-bundle over a
    closed surface $S$, as described in Definition~\ref{d-thinibundle}.
    We require this thin $I$-bundle to have no twists, so that it has
    two parallel boundary components each representing the same surface $S$.

    We then have the option of layering additional tetrahedra upon
    these boundary surfaces in order to change the curves formed by
    the boundary edges; see Definition~\ref{d-layering} for a detailed
    description of the layering process.  Finally we identify the two
    boundary surfaces according to some homeomorphism, forming an
    $S$-bundle over the circle.
\end{defn}

We proceed to describe in detail some particular classes of layered surface
bundles that feature in the $\leq 7$-tetrahedron non-orientable census.

\begin{defn} \label{d-surfacebundletypes}
    In Figure~\ref{fig-layeredsurfacebundles} we see five specific
    thin $I$-bundles.  Triangulations $T_6^1$, $T_6^2$ and $T_7$ are
    thin $I$-bundles over the torus and triangulations $K_6^1$ and $K_6^2$
    are thin $I$-bundles over the Klein bottle.  Each of these $I$-bundles
    is untwisted and therefore has two parallel boundary components,
    each formed from two triangles.

    \begin{figure}[htb]
    \psfrag{T61}{{\small $T_6^1$}} \psfrag{T62}{{\small $T_6^2$}}
    \psfrag{T7}{{\small $T_7$}}
    \psfrag{K61}{{\small $K_6^1$}} \psfrag{K62}{{\small $K_6^2$}}
    \psfrag{a1}{{\small $\alpha_1$}} \psfrag{a2}{{\small $\alpha_2$}}
    \psfrag{b1}{{\small $\beta_1$}} \psfrag{b2}{{\small $\beta_2$}}
    \centerline{\includegraphics[scale=0.7]{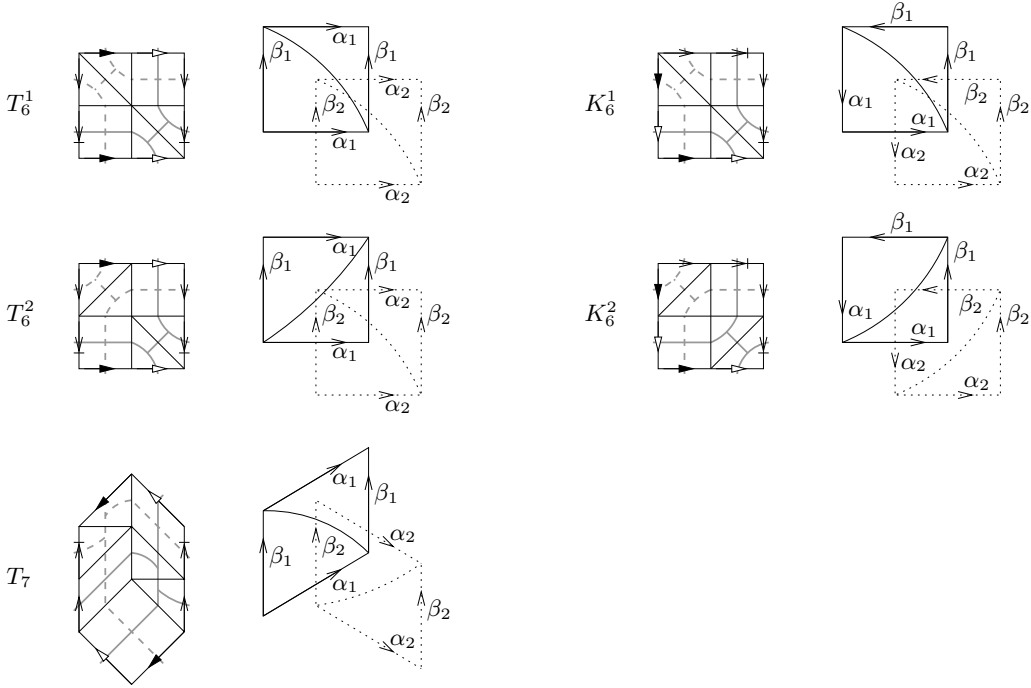}}
    \caption{The untwisted thin $I$-bundles
        $T_6^1$, $T_6^2$, $T_7$, $K_6^1$ and $K_6^2$}
    \label{fig-layeredsurfacebundles}
    \end{figure}

    The left hand portion of each diagram depicts the central
    surface decomposition of the thin $I$-bundle, complete with markings
    as described in Definition~\ref{d-markings}.  The right hand portion
    illustrates the upper and lower boundary components of the thin
    $I$-bundle, marked in solid lines and dotted lines respectively.
    Note that these boundary components correspond to the triangles of
    the central surface decomposition, since each tetrahedron enclosing
    a triangle provides a single boundary face.

    For each thin $I$-bundle we mark directed edges $\alpha_1$ and
    $\beta_1$ on the upper boundary component and $\alpha_2$ and
    $\beta_2$ on the lower boundary component.  Note that in each case
    $\alpha_1$ and $\beta_1$ generate the fundamental group of the upper
    boundary surface and $\alpha_2$ and $\beta_2$ generate
    the fundamental group of the lower boundary surface.

    Let $p$, $q$, $r$ and $s$ be integers and let $\theta$ be one of the thin
    $I$-bundles depicted in Figure~\ref{fig-layeredsurfacebundles}.
    We define $B_{\theta|p,q|r,s}$ to be the specific layered surface bundle
    obtained by identifying the upper and lower boundaries of $\theta$ so
    that directed edge $\alpha_1$ maps to $\alpha_2^p \beta_2^q$ and
    directed edge $\beta_1$ maps to $\alpha_2^r \beta_2^s$.

    Note that for some values of $p$, $q$, $r$ and $s$ this identification
    can be realised by an immediate mapping of the corresponding
    boundary faces.  On the other hand,
    for some values of $p$, $q$, $r$ and $s$ an additional layering of
    tetrahedra is required so that $\alpha_2^p \beta_2^q$,
    $\alpha_2^r \beta_2^s$ and the corresponding diagonal actually appear as
    edges of the lower boundary surface.

    Note also that for some values of $p$, $q$, $r$ and $s$ this
    construction is not possible since there is no homeomorphism
    identifying the upper and lower boundaries as required.
\end{defn}

\begin{example}
    To illustrate this construction we present the triangulations
    $B_{T_7|1,1|1,0}$ and $B_{T_6^2|-1,1|2,-1}$.

    The construction of $B_{T_7|1,1|1,0}$ begins with the thin $I$-bundle
    $T_7$, depicted in Figure~\ref{fig-t7-1110}.  Our task is to
    identify the upper and lower boundaries so that $\alpha_1$ and
    $\beta_1$ map to $\alpha_2 \beta_2$ and $\alpha_2$ respectively.

    \begin{figure}[htb]
    \psfrag{T61}{{\small $T_6^1$}} \psfrag{T62}{{\small $T_6^2$}}
    \psfrag{T7}{{\small $T_7$}}
    \psfrag{K61}{{\small $K_6^1$}} \psfrag{K62}{{\small $K_6^2$}}
    \psfrag{a1}{{\small $\alpha_1$}} \psfrag{a2}{{\small $\alpha_2$}}
    \psfrag{b1}{{\small $\beta_1$}} \psfrag{b2}{{\small $\beta_2$}}
    \psfrag{P}{{\small $P$}} \psfrag{Q}{{\small $Q$}}
    \psfrag{R}{{\small $R$}} \psfrag{S}{{\small $S$}}
    \psfrag{W}{{\small $W$}} \psfrag{X}{{\small $X$}}
    \psfrag{Y}{{\small $Y$}} \psfrag{Z}{{\small $Z$}}
    \centerline{\includegraphics[scale=0.7]{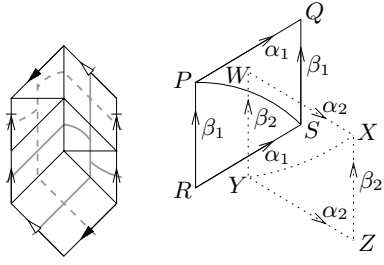}}
    \caption{Constructing the layered surface bundle $B_{T_7|1,1|1,0}$}
    \label{fig-t7-1110}
    \end{figure}

    This can in fact be done using a direct identification, by mapping
    boundary face {\em RPS} to {\em YZX} and mapping boundary face
    {\em PSQ} to {\em YWX}.  The final triangulation $B_{T_7|1,1|1,0}$
    therefore has seven tetrahedra, with no additional layering taking
    place.

    The construction of $B_{T_6^2|-1,1|2,-1}$ is slightly more complex.
    Figure~\ref{fig-t62-n112n1} illustrates the thin $I$-bundle $T_6^2$.
    Here we must identify the boundaries so that $\alpha_1$ and
    $\beta_1$ map to $\alpha_2^{-1}\beta_2$ and $\alpha_2^2\beta_2^{-1}$
    respectively.  Unfortunately this cannot be done using a direct
    identification of the boundaries since $\alpha_2^2\beta_2^{-1}$ does
    not appear as an edge of the lower boundary surface.

    \begin{figure}[htb]
    \psfrag{T61}{{\small $T_6^1$}} \psfrag{T62}{{\small $T_6^2$}}
    \psfrag{T7}{{\small $T_7$}}
    \psfrag{K61}{{\small $K_6^1$}} \psfrag{K62}{{\small $K_6^2$}}
    \psfrag{a1}{{\small $\alpha_1$}} \psfrag{a2}{{\small $\alpha_2$}}
    \psfrag{b1}{{\small $\beta_1$}} \psfrag{b2}{{\small $\beta_2$}}
    \psfrag{comb}{{\small $\alpha_2^{-1} \beta_2$}}
    \psfrag{P}{{\small $P$}} \psfrag{Q}{{\small $Q$}}
    \psfrag{R}{{\small $R$}} \psfrag{S}{{\small $S$}}
    \psfrag{W}{{\small $W$}} \psfrag{X}{{\small $X$}}
    \psfrag{Y}{{\small $Y$}} \psfrag{Z}{{\small $Z$}}
    \psfrag{A}{{\small $A$}} \psfrag{B}{{\small $B$}}
    \psfrag{C}{{\small $C$}} \psfrag{D}{{\small $D$}}
    \centerline{\includegraphics[scale=0.7]{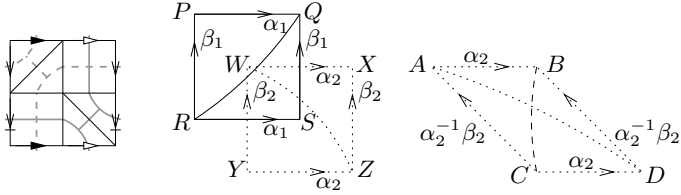}}
    \caption{Constructing the layered surface bundle $B_{T_6^2|-1,1|2,-1}$}
    \label{fig-t62-n112n1}
    \end{figure}

    We are therefore forced to layer a new tetrahedron onto the lower
    boundary.  This additional tetrahedron, labelled {\em ABCD} in the
    diagram, is layered upon edge {\em XZ} so that faces
    {\em YWZ} and {\em CBD} are identified and faces {\em WZX} and
    {\em ACB} are identified.  As a result we obtain a new lower boundary edge
    {\em AD} that indeed represents the curve $\alpha_2^2\beta_2^{-1}$.

    We can thus complete the triangulation by identifying the new lower
    boundary with the original upper boundary, mapping face
    {\em PQR} to {\em DBA} and face {\em QRS} to {\em DCA}.
    The final triangulation $B_{T_6^2|-1,1|2,-1}$ has seven tetrahedra.
\end{example}

For any layered surface bundle of a form described in
Definition~\ref{d-surfacebundletypes}, the underlying 3-manifold can be
identified using the following result.

\begin{theorem} \label{t-idsurfacebundle}
    For each set of integers $p$, $q$, $r$ and $s$ for which the
    corresponding triangulations can be constructed, the underlying
    3-manifolds of the layered surface bundles with parameters $p$, $q$,
    $r$ and $s$ are as follows.
    \begin{itemize}
        \item $B_{T_6^1|p,q|r,s}$ and $B_{T_6^2|p,q|r,s}$ are both
        triangulations of the space $\torus \times I / \homtwo{p}{r}{q}{s}$.

        \item $B_{T_7|p,q|r,s}$ is a triangulation of the space
        $\torus \times I / \homtwo{(p+q)}{(r+s)}{q}{s}$.

        \item Assume that $p+r$ is odd, $|p-q|=|r-s|=1$ and $p-q+r-s=0$.
        Then triangulations $B_{K_6^1|p,q|r,s}$ and $B_{K_6^2|p,q|r,s}$ both
        represent the space
        $\kb \times I / \homtwo{(\{p\}-\{r\})}{\{r\}}{0}{(s-r)}$,
        where the symbol $\{x\}$ is defined to be 1 if $x$ is odd and
        0 if $x$ is even.
    \end{itemize}
\end{theorem}

\begin{proof}
    It can be observed from Figure~\ref{fig-layeredsurfacebundles} that
    $\alpha_2$ is parallel to $\alpha_1$ and $\beta_2$ is parallel to
    $\beta_1$ within each of the $I$-bundles
    $T_6^1$, $T_6^2$, $K_6^1$ and $K_6^2$.  Within the
    $I$-bundle $T_7$ we find that $\alpha_2$ is parallel to $\alpha_1$
    and that $\beta_2$ is parallel to $\alpha_1 \beta_1$.

    Given these observations, it is a simple matter to convert the
    identification of the two boundary surfaces into the canonical form
    described in Section~\ref{s-defsurfacebundles} and thus establish
    the above results.
\end{proof}

\subsection{Plugged Thin $I$-Bundles} \label{s-pluggedthin}

Plugged thin $I$-bundles are formed by attaching layered solid
tori to twisted $I$-bundles over the torus.  The resulting 3-manifolds
are all Seifert fibred, where we allow {\sfslong}s over orbifolds as
well as over surfaces.  Details of the construction are as follows.

\begin{defn}[Plugged Thin $I$-Bundle] \label{d-pluggedthin}
    A {\em plugged thin $I$-bundle} is a 3-ma\-ni\-fold triangulation formed
    using the following construction.  Begin with one of the thin $I$-bundles
    over the torus depicted in Figure~\ref{fig-pluggedibundlestorus}.
    Note from the markings on the diagrams that each $I$-bundle is
    twisted and non-orientable, specifically with a twist as we wrap from
    top to bottom in each diagram and no twist as we wrap from left to right.

    \begin{figure}[htb]
    \psfrag{T61}{{\small $\twt_6^1$}}
    \psfrag{T62}{{\small $\twt_6^2$}}
    \psfrag{T63}{{\small $\twt_6^3$}}
    \psfrag{T64}{{\small $\twt_6^4$}}
    \centerline{\includegraphics[scale=0.7]{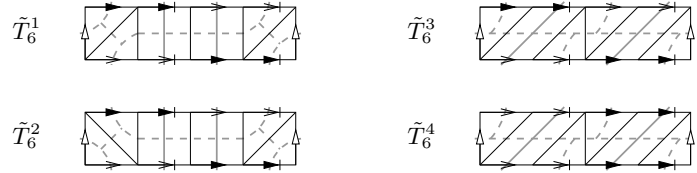}}
    \caption{The twisted thin $I$-bundles
        $\twt_6^1$, $\twt_6^2$, $\twt_6^3$ and $\twt_6^4$}
    \label{fig-pluggedibundlestorus}
    \end{figure}

    The four triangles of each central surface decomposition shown in
    Figure~\ref{fig-pluggedibundlestorus} correspond to the four
    boundary faces of each $I$-bundle.  In each case these boundary
    faces combine to form a torus as illustrated in
    Figure~\ref{fig-pluggedibundlestorusbdry}.  We observe that each
    of these torus boundaries is formed from two annuli, one on the left
    and one on the right.  Our construction is then completed by attaching a
    layered solid torus to each of these annuli as illustrated in
    Figure~\ref{fig-pluggedibundlestoruslst}.

    \begin{figure}[htb]
    \psfrag{T134}{{\small $\twt_6^1$, $\twt_6^3$, $\twt_6^4$}}
    \psfrag{T2}{{\small $\twt_6^2$}}
    \centerline{\includegraphics[scale=0.7]{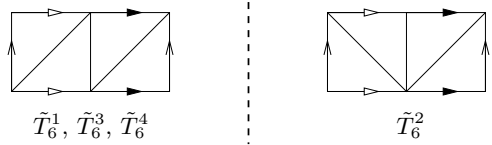}}
    \caption{The possible torus boundaries of thin $I$-bundles
        $\twt_6^1$, $\twt_6^2$, $\twt_6^3$ and $\twt_6^4$}
    \label{fig-pluggedibundlestorusbdry}
    \end{figure}

    \begin{figure}[htb]
    \psfrag{p1}{{\small $p_1$}}
    \psfrag{p2}{{\small $p_2$}}
    \psfrag{q1}{{\small $q_1$}}
    \psfrag{q2}{{\small $q_2$}}
    \psfrag{r1}{{\small $r_1$}}
    \psfrag{r2}{{\small $r_2$}}
    \centerline{\includegraphics[scale=0.7]{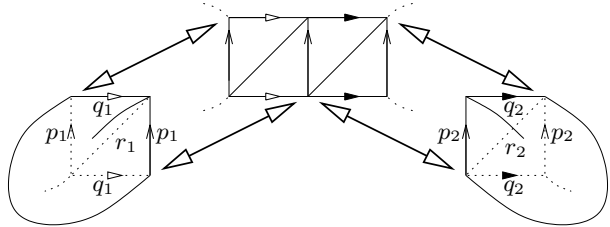}}
    \caption{Attaching layered solid tori to the torus boundary}
    \label{fig-pluggedibundlestoruslst}
    \end{figure}

    Let the layered solid tori have parameters $\lst(p_1,q_1,r_1)$ and
    $\lst(p_2,q_2,r_2)$ as explained in Definition~\ref{d-lstparams}.
    Furthermore, let the layered solid torus edges with parameters
    $p_i$ be attached to the left and right edges of the
    $I$-bundle boundary and let the layered solid torus edges with parameters
    $q_i$ be attached to the top and bottom edges of the
    $I$-bundle boundary as shown in
    Figure~\ref{fig-pluggedibundlestoruslst}.  Then the particular
    plugged thin $I$-bundle that has been constructed is denoted
    $H_{\theta|p_1,q_1|p_2,q_2}$, where $\theta$ denotes the original thin
    $I$-bundle chosen from Figure~\ref{fig-pluggedibundlestorus}.

    Note that instead of attaching a standard layered solid torus,
    the two faces of an annulus may simply be identified with
    each other by attaching the 0-tetrahedron degenerate
    $\lst(2,-1,-1)$, i.e.,
    a {\mobius} band.  For brevity, if a pair $p_i,q_i$ is omitted from
    the symbolic name of a plugged thin $I$-bundle then this pair is
    assumed to be $2,-1$.  For instance, the plugged thin $I$-bundle
    $H_{\twt_6^2|3,-1}$ is in fact $H_{\twt_6^2|3,-1|2,-1}$.
\end{defn}

Note that the triangulations $H_{\theta|p_1,q_1|p_2,q_2}$ and
$H_{\theta|p_2,q_2|p_1,q_1}$ are isomorphic.  This can be
seen from the symmetries of
the layered solid torus and of the thin $I$-bundles described in
Figure~\ref{fig-pluggedibundlestorus}.

The following result allows us to identify the underlying 3-manifold of
a plugged thin $I$-bundle.

\begin{theorem} \label{t-idpluggedibundletorus}
    Let $p_1$ and $q_1$ be coprime integers and let $p_2$ and $q_2$ be
    coprime integers, where $p_1 \neq 0$ and $p_2 \neq 0$.  Then the
    underlying 3-manifolds of the plugged thin $I$-bundles with
    parameters $p_1$, $q_1$, $p_2$ and $q_2$ are as
    follows.\footnote{In a previous version of this paper, the first 3-manifold
    was incorrectly given as $\sfs{\rpp}{(p_1,q_1),\ (p_2,q_2)}$.}

    \begin{itemize}
        \item $H_{\twt_6^1|p_1,q_1|p_2,q_2}$, $H_{\twt_6^2|p_1,q_1|p_2,q_2}$
        and $H_{\twt_6^3|p_1,q_1|p_2,q_2}$ are each triangulations of the
        {\sfslong} $\sfs{\rpp}{(p_1,q_1),\ (p_2,p_2 + q_2)}$.

        \item $H_{\twt_6^4|p_1,q_1|p_2,q_2}$ is a triangulation of the
        {\sfslong} $\sfs{\discref}{(p_1,q_1),\ (p_2,q_2)}$, where
        the orbifold $\discref$ is a disc with reflector boundary.
    \end{itemize}
\end{theorem}

\begin{proof}
    Consider the boundary torus of each of the thin $I$-bundles
    $\twt_6^1$, $\twt_6^2$, $\twt_6^3$ and $\twt_6^4$ as seen in
    Figure~\ref{fig-pluggedibundlestorusbdry}.
    We fill each of these boundary tori with
    circular fibres running parallel to the left and right
    sides as illustrated in Figure~\ref{fig-verticalbdrytorusfibres}.

    \begin{figure}[htb]
    \centerline{\includegraphics[scale=0.7]{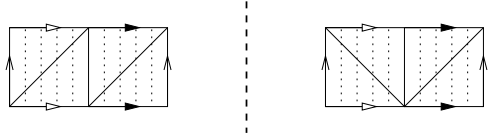}}
    \caption{Circular fibres on the boundary tori of thin $I$-bundles}
    \label{fig-verticalbdrytorusfibres}
    \end{figure}

    Our aim is to find compatible fibrations of the thin $I$-bundles.
    Each of the $I$-bundles $\twt_6^1$, $\twt_6^2$
    and $\twt_6^3$ can be realised as a trivial {\sfslong} over the
    {\mobius} band.
    The $I$-bundle $\twt_6^4$ on the other hand can be realised as a
    trivial {\sfslong} over an orbifold, where this base orbifold is an
    annulus with one reflector boundary.  In all cases the fibration of
    the $I$-bundle is compatible with the fibration of the boundary
    torus.

    Attaching our two layered solid tori then completes the
    fibrations, filling each base space with a disc
    and introducing exceptional fibres with parameters $(p_1,q_1)$ and
    $(p_2,p_2+q_2)$.  After normalising the Seifert invariants, the
    resulting 3-manifolds can be expressed as
    $\sfs{\rpp}{(p_1,q_1),\ (p_2,p_2+q_2)}$ and
    $\sfs{\discref}{(p_1,q_1),\ (p_2,q_2)}$ as claimed.
\end{proof}

\subsection{Plugged Thick $I$-Bundles} \label{s-pluggedthick}

Plugged thick $I$-bundles are similar to the plugged thin
$I$-bundles of Section~\ref{s-pluggedthin}, except that instead of
attaching layered solid tori directly to thin $I$-bundles we first wrap
the thin $I$-bundles with an additional padding of tetrahedra.
As with plugged thin $I$-bundles, the resulting 3-manifolds are all
Seifert fibred.

We begin by presenting four triangulations of a twisted $I$-bundle
over the torus, each of which has two vertices on the boundary.

\begin{defn}[Two-Vertex $I$-Bundles $\twt_5^1,\ldots,\twt_5^4$]
        \label{d-twovertexibundles}
    Let $\twt_3^1$, $\twt_3^2$ and $\twt_5^1$ be the thin $I$-bundles
    over the torus depicted in Figure~\ref{fig-smalltwistedtori}.
    Note that $\twt_3^1$ and $\twt_3^2$ are in fact the same
    triangulation presented in different ways.
    From the markings we see that each $I$-bundle is twisted and
    non-orientable, specifically with a twist as we wrap from top to bottom
    and no twist as we wrap from left to right.  The boundary torus of each
    $I$-bundle is illustrated alongside each diagram.

    \begin{figure}[htb]
    \psfrag{T51}{{\small $\twt_5^1$}}
    \psfrag{T31}{{\small $\twt_3^1$}}
    \psfrag{T32}{{\small $\twt_3^2$}}
    \centerline{\includegraphics[scale=0.7]{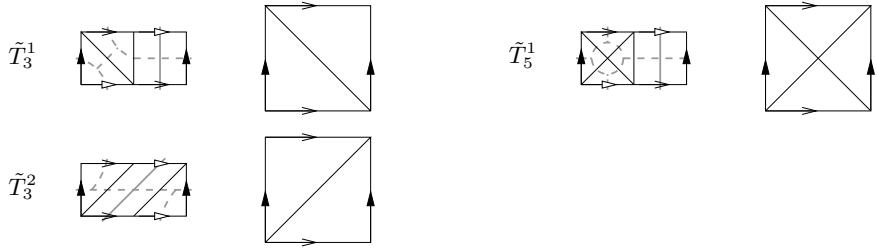}}
    \caption{The twisted thin $I$-bundles $\twt_3^1$, $\twt_3^2$ and
        $\twt_5^1$}
    \label{fig-smalltwistedtori}
    \end{figure}

    We see that $\twt_5^1$ already has two boundary vertices.  For
    $\twt_3^1$ and $\twt_3^2$ we modify the boundary by attaching a square
    pyramid formed from two tetrahedra.  The apex of this pyramid
    becomes the second boundary vertex.

    Figure~\ref{fig-pyramids} shows three new $I$-bundles $\twt_5^2$,
    $\twt_5^3$ and $\twt_5^4$ obtained in this fashion.  To construct
    $\twt_5^2$ and $\twt_5^3$ we attach a pyramid to the boundary of
    $\twt_3^1$; triangulation $\twt_5^3$ differs in that the base of
    the pyramid wraps around the upper and lower edges of the diagram.
    To construct $\twt_5^4$ we attach a pyramid to the boundary of $\twt_3^2$.
    The new two-vertex boundary tori are depicted on the right hand side of
    Figure~\ref{fig-pyramids}.

    \begin{figure}[htb]
    \psfrag{T52}{{\small $\twt_5^2$:}}
    \psfrag{T53}{{\small $\twt_5^3$:}}
    \psfrag{T54}{{\small $\twt_5^4$:}}
    \psfrag{T31}{{\small $\twt_3^1$}}
    \psfrag{T32}{{\small $\twt_3^2$}}
    \psfrag{+}{{\small $+$}}
    \centerline{\includegraphics[scale=0.7]{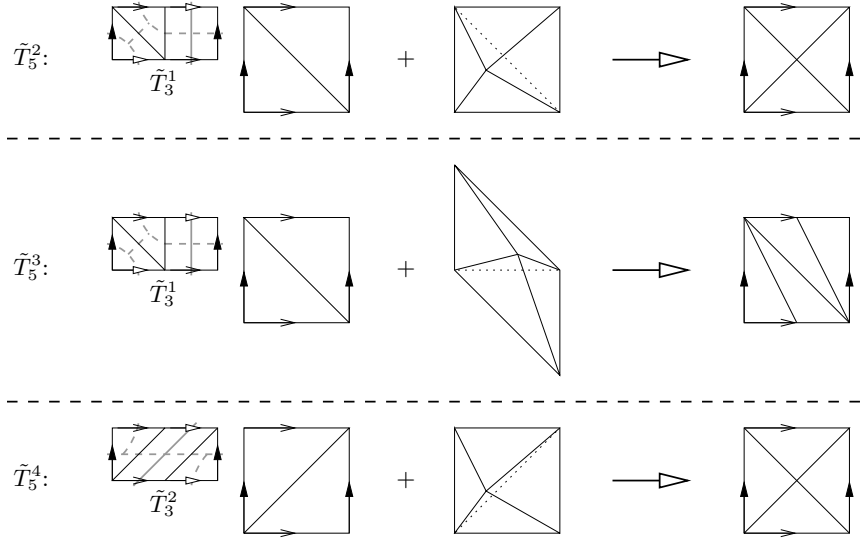}}
    \caption{Constructing twisted $I$-bundles $\twt_5^2$, $\twt_5^3$ and
        $\twt_5^4$}
    \label{fig-pyramids}
    \end{figure}

    We see then that each of the $I$-bundles $\twt_5^1$, $\twt_5^2$,
    $\twt_5^3$ and $\twt_5^4$ is formed from five tetrahedra and has a
    two-vertex torus boundary.
\end{defn}

Having constructed the $I$-bundles $\twt_5^1$, $\twt_5^2$, $\twt_5^3$ and
$\twt_5^4$, we proceed to define a plugged thick $I$-bundle as follows.

\begin{defn}[Plugged Thick $I$-Bundle] \label{d-pluggedthick}
    A {\em plugged thick $I$-bundle} is a 3-mani\-fold
    triangulation formed using the
    following construction.  Beginning with one of the
    two-vertex twisted $I$-bundles
    $\twt_5^1$, $\twt_5^2$, $\twt_5^3$ or $\twt_5^4$, we layer a single
    tetrahedron onto a specific edge of the boundary torus.
    This layering must form a new boundary edge running vertically from
    top to bottom;
    Figure~\ref{fig-layertwistedtori5} shows where this layering occurs
    for each of the $I$-bundles $\twt_5^1$, $\twt_5^2$, $\twt_5^3$ and
    $\twt_5^4$.

    \begin{figure}[htb]
    \psfrag{T5124}{{\small $\twt_5^1$, $\twt_5^2$, $\twt_5^4$}}
    \psfrag{T'5124}{{\small $\twt_5^{1'}$, $\twt_5^{2'}$, $\twt_5^{4'}$}}
    \psfrag{T53}{{\small $\twt_5^3$}} \psfrag{T'53}{{\small $\twt_5^{3'}$}}
    \centerline{\includegraphics[scale=0.7]{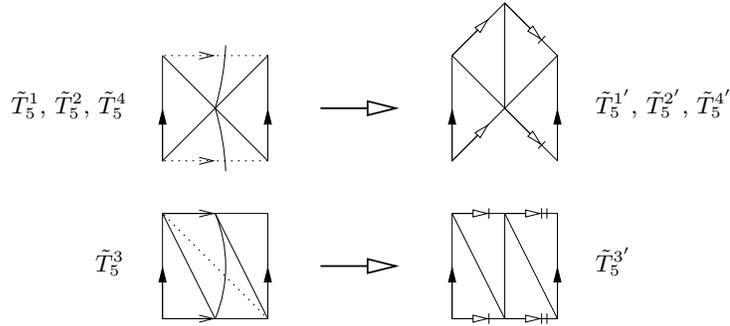}}
    \caption{Layering a tetrahedron to form twisted $I$-bundles
        $\twt_5^{1'}$, $\twt_5^{2'}$, $\twt_5^{3'}$ and $\twt_5^{4'}$}
    \label{fig-layertwistedtori5}
    \end{figure}

    Specifically, for $\twt_5^1$, $\twt_5^2$ and $\twt_5^4$
    the new tetrahedron is layered upon the top and bottom edges of the
    diagram.  For $\twt_5^3$ the new tetrahedron is layered upon the
    main diagonal.  The resulting $I$-bundles are labelled
    $\twt_5^{1'}$, $\twt_5^{2'}$, $\twt_5^{3'}$ and $\twt_5^{4'}$ respectively,
    and their new torus boundaries are shown on the right hand side of
    Figure~\ref{fig-layertwistedtori5}.

    At last we find ourselves in familiar territory.  The boundary tori
    of the new six-tetrahedron $I$-bundles $\twt_5^{1'}$, $\twt_5^{2'}$,
    $\twt_5^{3'}$ and $\twt_5^{4'}$ are each formed from two annuli, one on
    the left and one on the right.  As with the plugged thin $I$-bundles
    of the previous section, we complete the construction by
    attaching a layered solid torus to each annulus as illustrated in
    Figure~\ref{fig-pluggedthicklst}.

    \begin{figure}[htb]
    \psfrag{p1}{{\small $p_1$}}
    \psfrag{p2}{{\small $p_2$}}
    \psfrag{q1}{{\small $q_1$}}
    \psfrag{q2}{{\small $q_2$}}
    \psfrag{r1}{{\small $r_1$}}
    \psfrag{r2}{{\small $r_2$}}
    \centerline{\includegraphics[scale=0.7]{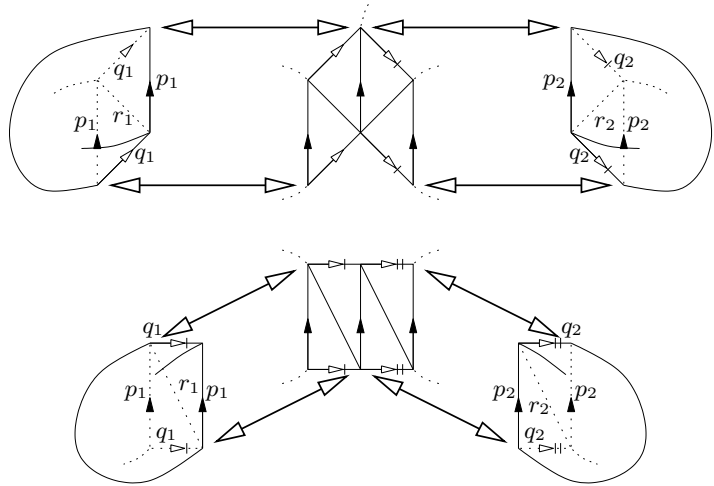}}
    \caption{Attaching layered solid tori to the torus boundary}
    \label{fig-pluggedthicklst}
    \end{figure}

    Let the layered solid tori have parameters $\lst(p_1,q_1,r_1)$ and
    $\lst(p_2,q_2,r_2)$, where the edges with parameters
    $p_i$ are attached to the left and right edges of the
    $I$-bundle boundary and where the edges with parameters
    $q_i$ are attached to the top and bottom edges of the
    $I$-bundle boundary.  Let $\theta$ be
    $\twt_5^1$, $\twt_5^2$, $\twt_5^3$ or $\twt_5^4$ according to which
    of the two-vertex twisted $I$-bundles was selected at the the
    beginning of the construction.
    Then the specific plugged thick $I$-bundle that has been constructed
    is denoted $K_{\theta|p_1,q_1|p_2,q_2}$.

    Again the two faces of an annulus may simply be identified with
    each other by attaching the degenerate $\lst(2,-1,-1)$, i.e.,
    a {\mobius} band.  If a pair $p_i,q_i$ is omitted from
    the symbolic name of a plugged thick $I$-bundle then this pair is
    once more assumed to be $2,-1$.  As an example, the plugged thick
    $I$-bundle $K_{\twt_5^3|3,-1}$ is in reality $K_{\twt_5^3|3,-1|2,-1}$.
\end{defn}

The identification of the underlying 3-manifold of a plugged thick
$I$-bundle is similar to that for a plugged thin $I$-bundle as seen
in the following result.

\begin{theorem} \label{t-idpluggedthick}
    Let $p_1$ and $q_1$ be coprime integers and let $p_2$ and $q_2$ be
    coprime integers, where $p_1 \neq 0$ and $p_2 \neq 0$.  Then the
    underlying 3-manifolds of the plugged thick $I$-bundles with
    parameters $p_1$, $q_1$, $p_2$ and $q_2$ are as
    follows.\footnote{In a previous version of this paper,
    these 3-manifolds were incorrectly given as
    $\sfs{\rpp}{(p_1,q_1),\ (p_2,q_2)}$ and
    $\sfs{\discref}{(p_1,q_1),\ (p_2,q_2)}$.}

    \begin{itemize}
        \item $K_{\twt_5^1|p_1,q_1|p_2,q_2}$, $K_{\twt_5^2|p_1,q_1|p_2,q_2}$
        and $K_{\twt_5^3|p_1,q_1|p_2,q_2}$ are each triangulations of the
        {\sfslong} $\sfs{\rpp}{(p_1,q_1),\ (p_2,p_2+q_2)}$.

        \item $K_{\twt_5^4|p_1,q_1|p_2,q_2}$ is a triangulation of the
        {\sfslong} $\sfs{\discref}{(p_1,q_1),\ (p_2,-q_2)}$, where
        the orbifold $\discref$ is a disc with reflector boundary.
    \end{itemize}
\end{theorem}

\begin{proof}
    The proof is almost identical to the proof of
    Theorem~\ref{t-idpluggedibundletorus}.
    Once more we fill the boundary tori of
    $\twt_5^{1'}$, $\twt_5^{2'}$, $\twt_5^{3'}$ and $\twt_5^{4'}$ with circular
    fibres that run parallel to the left and right sides
    as illustrated in Figure~\ref{fig-verticalbdrytorusfibresthick}.

    \begin{figure}[htb]
    \centerline{\includegraphics[scale=0.7]{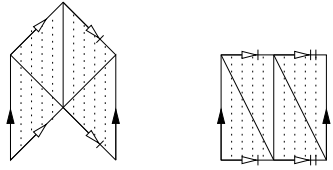}}
    \caption{Circular fibres on the boundary tori of thick $I$-bundles}
    \label{fig-verticalbdrytorusfibresthick}
    \end{figure}

    Compatible fibrations of the interior $I$-bundles and the
    exterior layered solid tori are found as before,
    resulting in the 3-manifolds listed above.
\end{proof}

\subsection{Exceptional Triangulations} \label{s-exceptional}

Three triangulations appear in the $\leq 7$-tetrahedron non-orientable
census that do not fit neatly into any of the families described thus far.
Each of these triangulations consists of precisely six tetrahedra and is
formed using a variant of a previous construction.  Specific details
of these three triangulations are as follows.

\begin{defn}[Triangulations $E_{6,1}$ and $E_{6,2}$]
    Triangulations $E_{6,1}$ and $E_{6,2}$ use a construction similar to
    the plugged thin $I$-bundles discussed in Section~\ref{s-pluggedthin}.
    Let $\twt_6^1$ and
    $\twt_6^2$ be the thin $I$-bundles over the torus depicted in
    the upper section of Figure~\ref{fig-e61-e62-base} (these are the
    same $\twt_6^1$ and $\twt_6^2$ as
    used in Definition~\ref{d-pluggedthin}).
    Each of these thin $I$-bundles has a torus boundary, illustrated
    in the lower portion of the figure.
    Triangulations $E_{6,1}$ and $E_{6,2}$ are formed from
    $\twt_6^1$ and $\twt_6^2$ respectively by identifying the faces of
    their torus boundaries as follows.

    \begin{figure}[htb]
    \psfrag{T61}{{\small $\twt_6^1$}}
    \psfrag{T62}{{\small $\twt_6^2$}}
    \psfrag{Bdry61}{{\small $\twt_6^1$ boundary}}
    \psfrag{Bdry62}{{\small $\twt_6^2$ boundary}}
    \centerline{\includegraphics[scale=0.7]{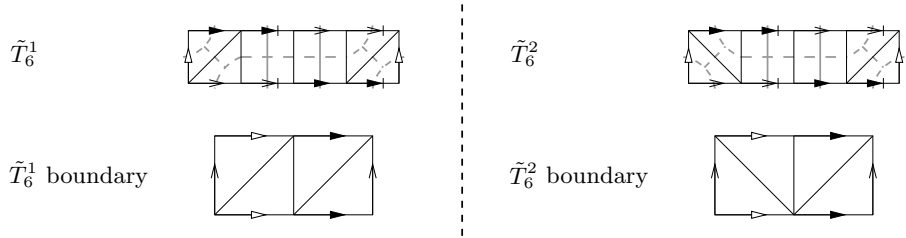}}
    \caption{Thin $I$-bundles used in the construction of
        $E_{6,1}$ and $E_{6,2}$}
    \label{fig-e61-e62-base}
    \end{figure}

    Figure~\ref{fig-e61} illustrates the construction of $E_{6,1}$.
    Beginning with $\twt_6^1$, we identify the boundary face
    {\em ADB} with {\em EBC} and we identify the boundary face
    {\em DBE} with {\em EFC}.  The resulting edge identifications are
    shown on the right hand side of the diagram.

    \begin{figure}[htb]
    \psfrag{A}{{\small $A$}} \psfrag{B}{{\small $B$}} \psfrag{C}{{\small $C$}}
    \psfrag{D}{{\small $D$}} \psfrag{E}{{\small $E$}} \psfrag{F}{{\small $F$}}
    \centerline{\includegraphics[scale=0.7]{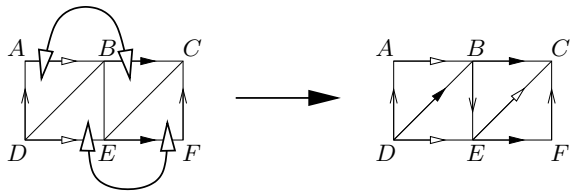}}
    \caption{Constructing exceptional triangulation $E_{6,1}$}
    \label{fig-e61}
    \end{figure}

    The construction of $E_{6,2}$ is illustrated in
    Figure~\ref{fig-e62}.  This time we begin with $\twt_6^2$,
    identifying the boundary face {\em AEB} with {\em ECF} and
    identifying the boundary face {\em DAE} with {\em BEC}.  Again the
    resulting edge identifications are shown.

    \begin{figure}[htb]
    \psfrag{A}{{\small $A$}} \psfrag{B}{{\small $B$}} \psfrag{C}{{\small $C$}}
    \psfrag{D}{{\small $D$}} \psfrag{E}{{\small $E$}} \psfrag{F}{{\small $F$}}
    \centerline{\includegraphics[scale=0.7]{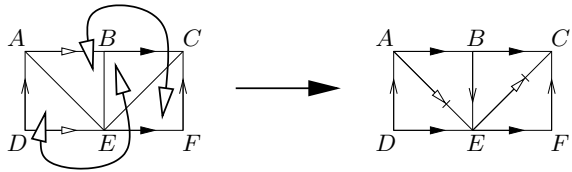}}
    \caption{Constructing exceptional triangulation $E_{6,2}$}
    \label{fig-e62}
    \end{figure}
\end{defn}

\begin{defn}[Triangulation $E_{6,3}$]
    Triangulation $E_{6,3}$ is formed from a pair of three-tet\-ra\-he\-dron
    thin $I$-bundles as follows.  Let $\twt_3^1$ be the thin $I$-bundle
    over the torus depicted in Figure~\ref{fig-e63-base} (this is the
    same $\twt_3^1$ as used in Definition~\ref{d-twovertexibundles}).
    This thin $I$-bundle has a torus boundary as shown in the right hand
    portion of the diagram.

    \begin{figure}[htb]
    \psfrag{T31}{{\small $\twt_3^1$}}
    \psfrag{Bdry31}{{\small $\twt_3^1$ boundary}}
    \centerline{\includegraphics[scale=0.7]{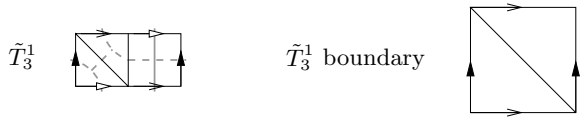}}
    \caption{The thin $I$-bundle used in the construction of $E_{6,3}$}
    \label{fig-e63-base}
    \end{figure}

    The construction of $E_{6,3}$ involves taking two copies of
    $\twt_3^1$ and identifying their boundary tori according to a
    particular homeomorphism.  This identification is illustrated in
    Figure~\ref{fig-e63}.  Specifically, face {\em ABD} is identified
    with {\em XZW} and face {\em ACD} is identified with {\em ZWY}.
    The resulting edge identifications are shown in the diagram.

    \begin{figure}[htb]
    \psfrag{A}{{\small $A$}} \psfrag{B}{{\small $B$}}
    \psfrag{C}{{\small $C$}} \psfrag{D}{{\small $D$}}
    \psfrag{W}{{\small $W$}} \psfrag{X}{{\small $X$}}
    \psfrag{Y}{{\small $Y$}} \psfrag{Z}{{\small $Z$}}
    \centerline{\includegraphics[scale=0.7]{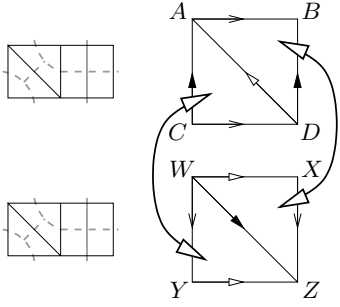}}
    \caption{Constructing exceptional triangulation $E_{6,3}$}
    \label{fig-e63}
    \end{figure}
\end{defn}

Each of the exceptional triangulations $E_{6,1}$, $E_{6,2}$ and $E_{6,3}$
can be converted using the elementary moves of
Section~\ref{s-algmanalysis-eltmoves} into a layered surface bundle, at
which point the underlying 3-manifold can be identified using
Theorem~\ref{t-idsurfacebundle}.  A list of
the resulting 3-manifolds can be found in Section~\ref{s-3manifolds}.

\section{Census Results} \label{s-census}

We conclude with a presentation of all closed non-orientable minimal
{\ppirr} triangulations formed from at most seven tetrahedra.
Recall from Section~\ref{s-intro} that this list contains 41 distinct
triangulations, together representing just eight different 3-manifolds.

Section~\ref{s-triangulations} lists these 41 triangulations according to
method of construction and Section~\ref{s-3manifolds} groups them
according to their underlying 3-manifolds.
Note that almost all of the 3-manifolds found in the census
allow for more than one minimal triangulation.  In such cases every
minimal triangulation is presented.

A {\regina} data file containing all of the census triangulations
listed here can be downloaded from the {\regina} website \cite{regina}.

\subsection{Triangulations} \label{s-triangulations}

We present here the 41 census triangulations ordered first by number of
tetrahedra and then by method of construction.
Table~\ref{tab-censusfrequencies} shows how these 41 triangulations
are distributed amongst the different families described in
Section~\ref{s-families}.  Note that the figures in the
six-tetrahedron column sum to 25 whereas the total is listed as 24; this
is because one of the six-tetrahedron triangulations can be viewed as
both a layered torus bundle and a layered Klein bottle bundle.

\renewcommand{\arraystretch}{1.2}
\begin{table}[htb]
\begin{center} \begin{tabular}{|l||r|r|r||r|}
    \hline
    \bf Tetrahedra & 1--5 & 6 & 7 & \bf Total \\
    \hline \hline
    Layered torus bundles & 0 & 6 & 4 & 10 \\
    Layered Klein bottle bundles & 0 & 8 & 0 & 8 \\
    Plugged thin $I$-bundles & 0 & 4 & 6 & 10 \\
    Plugged thick $I$-bundles & 0 & 4 & 7 & 11 \\
    Exceptional triangulations & \phantom{0}0 & 3 & 0 & 3 \\
    \hline \hline
    \bf Total & 0 & 24 & 17 & 41 \\
    \hline
\end{tabular} \end{center}
\caption{Frequencies of triangulations from different families}
\label{tab-censusfrequencies}
\end{table}

Since there are no closed non-orientable minimal {\ppirr} triangulations
formed from five tetrahedra or fewer, we use six tetrahedra as the
starting point for our detailed enumeration.

\subsubsection{Six Tetrahedra}

The six-tetrahedron closed non-orientable minimal {\ppirr}
triangulations are as follows.

\begin{itemize}
    \item The layered torus bundles
    $B_{T_6^1 | -1, 0 | -1, 1}$,
    $B_{T_6^1 |  0,-1 | -1, 0}$,
    $B_{T_6^1 |  0, 1 |  1, 0}$,
    $B_{T_6^1 |  1, 0 |  1,-1}$,
    $B_{T_6^2 | -1, 1 |  1, 0}$ and
    $B_{T_6^2 |  1, 0 |  0,-1}$
    as described by Definition~\ref{d-surfacebundletypes},
    where one of these layered torus bundles is isomorphic to a layered
    Klein bottle bundle as discussed below;

    \item The layered Klein bottle bundles
    $B_{K_6^1 | -1,0 | 0,-1}$,
    $B_{K_6^1 | 0,-1 | -1,0}$,
    $B_{K_6^1 | 0,1 | 1,0}$,
    $B_{K_6^1 | 1,0 | 0,1}$,
    $B_{K_6^2 | -1,0 | 0,-1}$,
    $B_{K_6^2 | 0,-1 | -1,0}$,
    $B_{K_6^2 | 0,1 | 1,0}$ and
    $B_{K_6^2 | 1,0 | 0,1}$
    as described by Definition~\ref{d-surfacebundletypes},
    where one of these layered Klein bottle bundles is isomorphic to a
    layered torus bundle as discussed below;

    \item The plugged thin $I$-bundles
    $H_{\twt_6^1}$,
    $H_{\twt_6^2}$,
    $H_{\twt_6^3}$ and
    $H_{\twt_6^4}$
    as described by Definition~\ref{d-pluggedthin};

    \item The plugged thick $I$-bundles
    $K_{\twt_5^1}$,
    $K_{\twt_5^2}$,
    $K_{\twt_5^3}$ and
    $K_{\twt_5^4}$
    as described by Definition~\ref{d-pluggedthick};

    \item The exceptional triangulations $E_{6,1}$, $E_{6,2}$ and
    $E_{6,3}$ as described in Section~\ref{s-exceptional}.
\end{itemize}

Although $25$ triangulations are named above, two of these are
isomorphic.  Specifically, the layered torus bundle
$B_{T_6^1 | 0,1 | 1,0}$ is isomorphic to the layered Klein bottle
bundle $B_{K_6^2 | 0,-1 | -1,0}$, leaving $24$ distinct triangulations in
our list.

\subsubsection{Seven Tetrahedra}

The seven-tetrahedron closed non-orientable minimal {\ppirr}
triangulations are as follows.

\begin{itemize}
    \item The layered torus bundles
    $B_{T_6^2 | -1,1 | 2,-1}$,
    $B_{T_6^2 | 0,-1 | -1,2}$,
    $B_{T_7 | -1,-1 | -1,0}$ and
    $B_{T_7 | 1,1 | 1,0}$
    as described by Definition~\ref{d-layeredsurfacebundle};

    \item The plugged thin $I$-bundles
    $H_{\twt_6^1 | 3,-2}$,
    $H_{\twt_6^1 | 3,-1}$,
    $H_{\twt_6^2 | 3,-2}$,
    $H_{\twt_6^2 | 3,-1}$,
    $H_{\twt_6^3 | 3,-1}$ and
    $H_{\twt_6^4 | 3,-1}$
    as described by Definition~\ref{d-pluggedthin};

    \item The plugged thick $I$-bundles
    $K_{\twt_5^1 | 3,-1}$,
    $K_{\twt_5^2 | 3,-2}$,
    $K_{\twt_5^2 | 3,-1}$,
    $K_{\twt_5^3 | 3,-2}$,
    $K_{\twt_5^3 | 3,-1}$,
    $K_{\twt_5^4 | 3,-2}$ and
    $K_{\twt_5^4 | 3,-1}$
    as described by Definition~\ref{d-pluggedthick}.
\end{itemize}

\subsection{3-Manifolds} \label{s-3manifolds}

We close with a table of all closed non-orientable {\ppirr}
3-manifolds formed from seven tetrahedra or fewer.
These 3-manifolds are listed in Table~\ref{tab-3mfds}, along with their
first homology groups and minimal triangulations.  Recall from
Section~\ref{s-families} that the orbifold $\discref$ is a disc with
reflector boundary.

\renewcommand{\arraystretch}{1.5}
\begin{table}[htb]
\small
\begin{center} \begin{tabular}{|c|l|l|l|}
    \hline
    $\Delta$ & \bf 3-Manifold & \bf Triangulations & \bf Homology \\
    \hline \hline
6 & $\torus \times I / \homtwo{-1}{1}{1}{0}$ &
    $B_{T_6^2 | -1,1|1,0}$
    & $\Z$ \\

  \cline{2-4}
  & $\torus \times I / \homtwo{0}{1}{1}{0}$ &
    $B_{T_6^1 | -1,0|-1,1}$,
    $B_{T_6^1 | 0,-1|-1,0}$,
    & $\Z \oplus \Z$ \\
  & &
    $(B_{T_6^1 | 0,1|1,0} = B_{K_6^2 | 0,-1|-1,0})$,
    & \\
  & &
    $B_{T_6^1 | 1,0|1,-1}$,
    $B_{K_6^1 | 0,-1|-1,0}$,
    & \\
  & &
    $E_{6,3}$
    & \\

  \cline{2-4}
  & $\kb \times \scircle$ &
    $B_{T_6^2 | 1,0|0,-1}$,
    $B_{K_6^1 | 1,0|0,1}$,
    & $\Z \oplus \Z \oplus \Z_2$ \\
  & &
    $B_{K_6^2 | 1,0|0,1}$
    & \\

  \cline{2-4}
  & $\kb \times I / \homtwo{-1}{1}{0}{-1}$ &
    $B_{K_6^1 | 0,1|1,0}$,
    $B_{K_6^2 | 0,1|1,0}$,
    & $\Z \oplus \Z_4$ \\
  & &
    $H_{\twt_6^1}$,
    $H_{\twt_6^2}$,
    $H_{\twt_6^3}$,
    & \\
  & &
    $K_{\twt_5^1}$,
    $K_{\twt_5^2}$,
    $K_{\twt_5^3}$,
    $E_{6,2}$
    & \\

  \cline{2-4}
  & $\kb \times I / \homtwo{1}{0}{0}{-1}$ &
    $B_{K_6^1 | -1,0|0,-1}$,
    $B_{K_6^2 | -1,0|0,-1}$,
    & $\Z \oplus \Z_2 \oplus \Z_2$ \\
  & &
    $H_{\twt_6^4}$,
    $K_{\twt_5^4}$,
    $E_{6,1}$
    & \\
    \hline
    \hline

7 & $\torus \times I / \homtwo{2}{1}{1}{0}$ &
    $B_{T_6^2 | -1,1 | 2,-1}$,
    $B_{T_6^2 | 0,-1 | -1,2}$,
    & $\Z \oplus \Z_2$ \\
  & &
    $B_{T_7 | -1,-1 | -1,0}$,
    $B_{T_7 | 1,1 | 1,0}$
    & \\
  \cline{2-4}
  & $\sfs{\rpp}{(2,1)\ (3,1)}$ &
    $H_{\twt_6^1 | 3,-2}$,
    $H_{\twt_6^1 | 3,-1}$,
    $H_{\twt_6^2 | 3,-2}$,
    & $\Z$ \\
  & &
    $H_{\twt_6^2 | 3,-1}$,
    $H_{\twt_6^3 | 3,-1}$,
    $K_{\twt_5^1 | 3,-1}$,
    & \\
  & &
    $K_{\twt_5^2 | 3,-2}$,
    $K_{\twt_5^2 | 3,-1}$,
    $K_{\twt_5^3 | 3,-2}$,
    & \\
  & &
    $K_{\twt_5^3 | 3,-1}$
    & \\
  \cline{2-4}
  & $\sfs{\discref}{(2,1)\ (3,1)}$ &
    $H_{\twt_6^4 | 3,-1}$,
    $K_{\twt_5^4 | 3,-2}$,
    $K_{\twt_5^4 | 3,-1}$
    & $\Z \oplus \Z_2$ \\
\hline
\end{tabular} \end{center}
\caption{All eight 3-manifolds and their 41 minimal triangulations}
\label{tab-3mfds}
\end{table}

\bibliographystyle{amsplain}
\bibliography{pure}

\newcommand{\noopsort}[1]{}
\providecommand{\bysame}{\leavevmode\hbox to3em{\hrulefill}\thinspace}
\providecommand{\MR}{\relax\ifhmode\unskip\space\fi MR }
\providecommand{\MRhref}[2]{%
  \href{http://www.ams.org/mathscinet-getitem?mr=#1}{#2}
}
\providecommand{\href}[2]{#2}
\begin{thebibliography}{10}

\bibitem{italian-nor6}
Gennaro Amendola and Bruno Martelli, \emph{Non-orientable 3-manifolds of small
  complexity}, Topology Appl. \textbf{133} (2003), no.~2, 157--178.

\bibitem{italian-nor7}
\bysame, \emph{Non-orientable 3-manifolds of complexity up to 7}, Topology
  Appl. \textbf{150} (2005), no.~1-3, 179--195.

\bibitem{regina}
Benjamin~A. Burton, \emph{Regina: Normal surface and 3-manifold topology
  software}, {\tt http://\allowbreak regina.\allowbreak sourceforge.\allowbreak
  net/}, 1999--2005.

\bibitem{burton-thesis}
\bysame, \emph{Minimal triangulations and normal surfaces}, Ph.D. thesis,
  University of Melbourne, 2003, available from {\tt http://\allowbreak
  regina.\allowbreak sourceforge.\allowbreak net/}.

\bibitem{burton-facegraphs}
\bysame, \emph{Face pairing graphs and 3-manifold enumeration}, J. Knot Theory
  Ramifications \textbf{13} (2004), no.~8, 1057--1101.

\bibitem{burton-regina}
\bysame, \emph{Introducing {R}egina, the 3-manifold topology software},
  Experiment. Math. \textbf{13} (2004), no.~3, 267--272.

\bibitem{cuspedcensus}
Patrick~J. Callahan, Martin~V. Hildebrand, and Jeffrey~R. Weeks, \emph{A census
  of cusped hyperbolic 3-manifolds}, Math. Comp. \textbf{68} (1999), no.~225,
  321--332.

\bibitem{hemion}
Geoffrey Hemion, \emph{The classification of knots and 3-dimensional spaces},
  Oxford Science Publications, Oxford University Press, Oxford, 1992.

\bibitem{hempel}
John Hempel, \emph{3-manifolds}, Annals of Mathematics Studies, no.~86,
  Princeton University Press, Princeton, NJ, 1976.

\bibitem{cuspedcensusold}
Martin~V. Hildebrand and Jeffrey~R. Weeks, \emph{A computer generated census of
  cusped hyperbolic 3-manifolds}, Computers and Mathematics (Cambridge, MA,
  1989), Springer, New York, 1989, pp.~53--59.

\bibitem{jacorubin-algorithms}
William Jaco, David Letscher, and J.~Hyam Rubinstein, \emph{Algorithms for
  essential surfaces in 3-manifolds}, Topology and Geometry: Commemorating
  SISTAG, Contemporary Mathematics, no. 314, Amer. Math. Soc., Providence, RI,
  2002, pp.~107--124.

\bibitem{0-efficiency}
William Jaco and J.~Hyam Rubinstein, \emph{0-efficient triangulations of
  3-manifolds}, J. Differential Geom. \textbf{65} (2003), no.~1, 61--168.

\bibitem{layeredlensspaces}
\bysame, \emph{Layered triangulations of lens spaces}, In preparation, 2003.

\bibitem{jacotollefson-algorithms}
William Jaco and Jeffrey~L. Tollefson, \emph{Algorithms for the complete
  decomposition of a closed {$3$}-manifold}, Illinois J. Math. \textbf{39}
  (1995), no.~3, 358--406.

\bibitem{kneser-normal}
Hellmuth Kneser, \emph{Geschlossene {F}l\"achen in dreidimensionalen
  {M}annigfaltigkeiten}, Jahresbericht der Deut. Math. Verein. \textbf{38}
  (1929), 248--260.

\bibitem{italian10}
Bruno Martelli, \emph{Complexity of 3-manifolds}, Preprint, {\tt
  math.\allowbreak GT/\allowbreak 0405250}, January 2005.

\bibitem{italian9}
Bruno Martelli and Carlo Petronio, \emph{Three-manifolds having complexity at
  most 9}, Experiment. Math. \textbf{10} (2001), no.~2, 207--236.

\bibitem{italianfamilies}
\bysame, \emph{Complexity of geometric three-manifolds}, Geom. Dedicata
  \textbf{108} (2004), no.~1, 15--69.

\bibitem{matveev6}
Sergei~V. Matveev, \emph{Tables of 3-manifolds up to complexity 6},
  Max-Planck-Institut f\"{u}r Mathematik Preprint Series (1998), no.~67,
  available from {\tt http://www.\allowbreak mpim-bonn.\allowbreak
  mpg.\allowbreak de/\allowbreak html/\allowbreak pre\-prints/\allowbreak
  preprints.html}.

\bibitem{matveev9}
\bysame, \emph{Computer classification of 3-manifolds}, Russ. J. Math. Phys.
  \textbf{7} (2000), no.~3, 319--329.

\bibitem{pachner-moves}
Udo Pachner, \emph{P.{L}. homeomorphic manifolds are equivalent by elementary
  shellings}, European J. Combin. \textbf{12} (1991), no.~2, 129--145.

\bibitem{schubert-normal}
Horst Schubert, \emph{Bestimmung der {P}rimfaktorzerlegung von {V}erkettungen},
  Math. Z. \textbf{76} (1961), 116--148.

\bibitem{turaev-viro}
Vladimir~G. Turaev and Oleg~Y. Viro, \emph{State sum invariants of
  {$3$}-manifolds and quantum {$6j$}-symbols}, Topology \textbf{31} (1992),
  no.~4, 865--902.

\end{thebibliography}

\end{document}